\documentclass[12pt]{article}
\usepackage{amsmath, amssymb, amsthm, url, hyperref}
\usepackage{geometry}
\title{Semi-slant and hemi-slant Riemannian Maps to K\"ahler Manifolds}
\author{}
\date{}
\newtheorem{theorem}{Theorem}[section]
\newtheorem{lemma}{Lemma}[section]

\newtheorem{definition}{Definition}[section]
\newtheorem{example}{Example}[section]
\theoremstyle{remark}

\date{}
\title{\textbf{Clairaut semi-slant/hemi-slant Riemannian maps to K\"ahler manifolds}}
\author{Jyoti Yadav, Gauree Shanker\thanks{corresponding author, Email: gauree.shanker@cup.edu.in}}
\begin{document}
\maketitle
\begin{abstract}
The aim of this paper is to study Clairaut semi-slant(hemi-slant) Riemannian maps to K\"ahler manifolds. We find the geodesic behavior of a regular curve on a base manifold of these maps. A central part of this study is devoted to investigate the necessary and sufficient conditions for these maps to be Clairaut semi-slant(hemi-slant) Riemannian maps. In this context, we discuss the necessary and sufficient conditions for these maps to be totally geodesic harmonic maps and examine the geometry of distributions. We also study Casorati inequality and Chen’s first inequality for semi-slant and hemi-slant Riemannian maps and discuss the equality case. To support the theory, some examples are provided.
\end{abstract}

\section{Introduction}
The concept of Riemannian maps, introduced by Fischer in 1992, provides a unifying framework that encompasses both Riemannian submersions and isometric immersions between Riemannian manifolds. These maps serve as an essential tool for analyzing and comparing the geometric structures of different manifolds. A notable feature of Riemannian maps is their connection with the generalized Eikonal equation, linking them to physical theories such as geometric optics and Maxwell–Schrödinger systems. This perspective highlights the relevance of Riemannian maps not only in differential geometry but also in mathematical physics and Lagrangian field theory. For further study, we refer  \cite{Falcitelli, Fischer, SahinBook, SahinRmaps}.\\
In this direction, the notion of Clairaut submersion, inspired by the classical Clairaut relation in differential geometry, was introduced by Bishop \cite{Bishop} and has since been studied in connection with various maps and manifolds. Researchers have generalized the Clairaut Theorem with various kind of maps on manifolds \cite{Meena_Zawadzki, Roy, Tastan, Maini, Gunduzalp, YadavJ, KumarCSIRm, YadavCIRM, SahinCRMap, YadavCAIRM, Polat, MeenaCRMap, MeenaCCRM, Shankeryadav, Yadavslant}. Motivated by these developments, in this paper, we investigate Clairaut semi-slant/hemi-slant Riemannian maps to K\"ahler manifolds. An important direction in submanifold theory is the study of relationships between intrinsic and extrinsic invariants. In this context, Chen (1993) established a fundamental inequality, widely known as Chen’s first inequality which relates the geometry of a submanifold to its ambient space. This inequality has become a powerful tool with applications in the study of geometric flows, rigidity theorems, curvature bounds, and classification of submanifolds. Several researchers \cite{sahinchen, MeenaChen, chen1} have extended Chen’s inequality to the setting of Riemannian maps to real, complex, K\"ahler, Sasakian, Kenmotsu, cosymplectic, generalized complex and generalized Sasakian space forms, motivating further exploration in broader geometric frameworks.\\
Parallel to this, the notion of Casorati curvature, originally introduced for surfaces in Euclidean three-space, has been generalized to submanifolds in Riemannian geometry as an extrinsic invariant. It is defined by the normalized squared length of the second fundamental form. It plays a central role in optimal inequalities, and recent studies have investigated such inequalities for various maps and ambient geometries \cite{Lee, cas cur}. Therefore, it is interesting to obtain optimal inequalities for the Casorati curvatures of submanifolds in different ambient spaces. The paper is organized as follows. In Section 2, we recall the necessary preliminaries. In section 3, we define Clairaut semi-slant Riemannian maps to K\"ahler manifolds. Further, we obtain the necessary and sufficient conditions for semi-invariant Riemannian maps to be Clairaut semi-slant Riemannian maps with the help of geodesic. Also, we address the harmonicity and totally geodesic conditions of these maps. 
Section 4 deals with Clairaut hemi-slant Riemannian maps to K\"ahler manifolds. Next, first we find the geodesic behaviour of a regular curve and obtain a necessary and sufficient condition for hemi-slant Riemannian maps to Clairaut hemi-slant Riemannian maps. Finally, we find a necessary and sufficient condition for these maps to be totally geodesic and investigate the geometry of distribution. In section 5, we discuss some curvature inequalities(Casorati curvature, Chen's first inequality) for semi-slant and hemi-slant Riemannian map. Finally, we provide conclusion and future scope of the current research work.
\section{Preliminaries}
In this section, we discuss some basic definitions and results required for this paper.\\
Let $F : (M^{m}, g_{M}) \to (N^{n}, g_{N})$ be a smooth map between two Riemannian manifolds, where $0 < \text{rank}F\leq \min \{ m, n \}$. For given $p \in M$, $\mathcal{V}_{p} = \ker F_{* p}$ denotes the vertical space and orthogonal complement of $\mathcal{V}_{p}$ is $\mathcal{H}_{p} = (\ker F_{* p})^{\perp}.$ Therefore, the tangent space $T_{p} M$ can be written as \cite{Fischer}
\[
T_{p} M = \ker F_{* p} \oplus (\ker F_{* p})^{\perp}.
\]
Moreover, range of $F_{* p}$ is denoted by $\mathrm{range} \, F_{* p}$ and its orthogonal complement is given by $(\mathrm{range} \, F_{* p})^{\perp}$. Thus tangent space at $T_{F(p)} N,~F(p) \in N$,  is the direct sum of $\mathrm{range} \, F_{* p}$ and $(\mathrm{range} \, F_{* p})^{\perp}$, that is,
\[
T_{F(p)} N = \mathrm{range} \, F_{* p} \oplus (\mathrm{range} \, F_{* p})^{\perp}.
\]

\begin{definition}\cite{Fischer}
A smooth map $F$ is said to be a \emph{Riemannian map} at $p \in M$ if the horizontal restriction 
\[
F_{* p} : (\ker F_{* p})^{\perp} \to \mathrm{range} \, F_{* p}
\]
is a linear isometry between the inner product spaces $(\ker F_{* p})^{\perp}, g_{M}|_{(\ker F_{* p})^{\perp}})$ and\\ $(\mathrm{range} \, F_{* p}, g_{N}|_{\mathrm{range} \, F_{* p}})$, that is,
\[
g_{N}(F_{*} X, F_{*} Y) = g_{M}(X, Y).
\]
\end{definition} 
For a Riemannian map $F,$ define $S_V$ as \cite{SahinBook}
\begin{equation}\label{Sv}
\nabla^N_{F_*X}V = -S_{V}F_*X+\nabla^{F\perp}_X V
\end{equation}
where $S_VF_*X$ is the tangential component (a vector field along $F$) of $\nabla^N_{F_*X}V.$ Thus, at $p\in M,$ we have $\nabla^N_{F_*X}V(p)\in T_{F(p)}N, S_VF_*X\in F_*(T_pM), \nabla^{F\perp}_X V(p)\in (F_*(T_pM))^\perp.$\\
Suppose $F:(M, g_M)\rightarrow(N, g_N)$ be a smooth map between Riemannian manifolds. Then the differential $F_*$ of $F$ can be viewed as a section of bundle \text{Hom}$(TM, F^{-1}TN)\rightarrow M,$ where $F^{-1}TN$ is the pullback bundle whose fibers at $p \in M$ is $(F^{-1}TN)_{p} =T_{F (p)}N, p \in M.$ The bundle \text{Hom}$(TM, F^{-1}TN)$ has a connection $\nabla$ induced from the Levi-Civita connection $\nabla^M$ and the pullback connection $\nabla^{N_{F}}.$ Then the second fundamental form of $F$ is given by \cite{Nore}
\begin{equation}\label{SFF}
(\nabla F_*)(X, Y) = \nabla^{N_{F}}_X F_*Y - F_*(\nabla^M_X Y)
\end{equation}
for all $X, Y \in\Gamma(TM).$
\begin{lemma}\cite{SahinBook}\label{URM}
Let $F:(M, g_M)\rightarrow (N, g_N)$ be a Riemannian map. Then $F$ is an umbilical Riemannian map if and only if
\begin{equation}
(\nabla F_*)(X, Y) = g_M(X, Y)H
\end{equation}
for $X, Y\in\Gamma(kerF_*)^\perp$ and $H$ is a vector field on $(rangeF_*)^\perp.$
\end{lemma}
\begin{definition}\cite{Yano}
Let $(M, g_M)$ be an almost Hermitian manifold, then $M$ admits a tensor $J$ of
type $(1, 1)$ such that $J^2 = -I$ and \begin{equation}\label{JXJY}
  g_M(JX, JY) = g_M(X, Y)  
\end{equation} for all $X, Y \in\Gamma(TM).$ An almost Hermitian manifold $M$ is called a K\"ahler manifold if 
\begin{equation}\label{nablaJ}
 (\nabla^M_X J)Y = 0   
\end{equation}for all $X, Y \in\Gamma(TM),$ where $\nabla^M$ is a Levi–Civita connection on $M.$
\end{definition}
\begin{definition}\cite{park}
A Riemannian map $F : (M, g_M)\rightarrow (N, g_N, J)$ from a Riemannian manifold to an almost Hermitian manifold is called a semi-slant Riemannian map if there is a distribution $D_1 \subset (kerF_*)^\perp$ such that
$$(kerF_*)^\perp = D_1\oplus D_2, J(F_*D_1)=F_*D_1,$$ and the angle $\theta = \theta(X)$ between $JF_*X$ and the space $F_*(D_2)_p$ is constant for non-zero vector $X\in(D_2)_p,~~p\in M,$ where $D_2$ is the orthogonal complement of $D_1$ in $(kerF_*)^\perp.$ The angle $\theta$ is known as semi-slant angle.
\end{definition}
Let $F : (M, g_M)\rightarrow(N, g_N, J)$ be a semi-slant Riemannian map. Then for
$F_*X\in\Gamma(range F_*),$ we have
\begin{equation}\label{SJU}
JF_*X = \phi F_*X+\omega F_*X,   
\end{equation}
where $\phi F_*X\in\Gamma(rangeF_*)$ and $\omega F_*X\in\Gamma(rangeF_*)^\perp$\\
and for $V\in\Gamma(rangeF_*)^\perp,$ we have
\begin{equation}\label{SJV}
JV = BV+CV,   
\end{equation}
where $BV\in\Gamma(rangeF_*)$ and $CV\in\Gamma(rangeF_*)^\perp.$
\begin{theorem}\cite{park}\label{phiFx}
Let $F$ be a semi-slant Riemannian map from a Riemannian manifold $(M, g_M)$ to an almost Hermitian manifold $(N, g_N, J)$ with the semi-slant angle $\theta.$
Then we obtain
$$\phi^2F_*X = -\cos^2\theta F_*X~~ for~~ X\in\Gamma(D_2).$$
\end{theorem}
\begin{definition}\cite{sahinhemi-slant}
A Riemannian map $F:(M, g_M)\rightarrow (N, g_N, J)$ between a Riemannian manifold to an almost Hermitian manifold is called a hemi-slant Riemannian map if $(rangeF_*)$ of $F$ admits two orthogonal complementary distributions $D^{\psi}$ and $D^{\perp}$ such that $D^{\psi}$ is slant and $D^{\perp}$ is an anti-invariant, that is, $$range F_* = D^{\psi}\oplus D^{\perp}.$$ In this case, the angle $\psi$ is called the hemi-slant angle of the Riemannian map.
\end{definition}
Let $F$ be a hemi-slant Riemannian map, then for $F_*X\in\Gamma (range F_*),$ we have
\begin{equation}\label{hemiJV}
 JF_*X = \phi F_*X+ \omega F_*X,   
\end{equation}
where $\phi F_*X\in\Gamma(rangeF_*)$ and $\omega F_*X\in\Gamma(rangeF_*)^\perp.$ Also, for any $V\in\Gamma(rangeF_*)^\perp,$ we have
\begin{equation}\label{hemiJxi}
 JV = BV+ CV,   
\end{equation} where $BV\in\Gamma(rangeF_*)$ and $CV\in\Gamma(rangeF_*)^\perp.$
\begin{theorem}\label{hemiangle}\cite{sahinhemi-slant}
Let $F$ be a Riemannian map from a Riemannian manifold $(M, g_M)$ to an almost Hermitian
manifold $(N, g_N, J).$ Then $F$ is a hemi-slant Riemannian map if
and only if there exists a constant $\lambda\in[-1, 0]$ and a distribution $D$ on \text{range}$F_*$ such that\\
\begin{itemize}
    \item [(a)] $D = \{F_*X\in\Gamma(rangeF_*)|\phi^2 F_*X = \lambda F_*X\},$
    \item [(b)] for any $F_*X\in\Gamma(rangeF_*)$ orthogonal to $D,$ we have $\phi F_*X = 0.$
\end{itemize}
Moreover, in this case $\lambda = -\cos^2\theta,$ where $\theta$ is the slant angle of $F.$
\end{theorem}
\begin{definition}\label{CRMps}\cite{MeenaCRMap}
	A Riemannian map $F : (M, g_M) \rightarrow (N, g_N)$ between Riemannian manifolds is a called Clairaut Riemannian map if there is a function $\tilde{s}: N \rightarrow \mathbb{R}^+$ such that for every geodesic $\beta$ on $N,$ the function 
	$(\tilde{s} \circ \beta) \sin \psi(t)$  is constant, where $\psi(t)$ is the angle between $\dot\beta(t)$ and horizontal subspace at $\beta(t).$
\end{definition}
\begin{theorem}\cite{MeenaCRMap}\label{NSC}
	Let $F : (M, g_M) \rightarrow (N, g_N)$ be a Riemannian map between Riemannian manifolds such that $(rangeF_*)^\perp$ is totally geodesic and $range F_*$ is connected, and let $\alpha, \beta = F\circ \alpha$ be geodesics on $M$ and $N,$ respectively.Then, $F$ is a Clairaut Riemannian map with $\tilde{s} = e^g$ if and only if one of the following conditions holds: 
\begin{itemize}
	\item[(i)] $S_V F_*X = -V(g) F_*X$, where $F_*X \in\Gamma (rangeF_*)$ and 
	$V \in\Gamma(rangeF_*)^\perp$ are components of  $\dot{\beta}(t)$. 
	\item [(ii)] $F$ is umbilical map, and has $H = -\nabla^N g$ , where $g$ is a smooth function on $N$ and $H$ is the mean curvature vector field of $rangeF_*$.
\end{itemize} 
\end{theorem}
For all $U, V\in\Gamma(rangeF_*)^\perp,$ we have \cite{Yadavbase}
$$\nabla^N_U V = R(\nabla^N_U V)+\nabla^{F\perp}_U V,$$ where $R(\nabla^N_U V)$ and $\nabla^{F\perp}_U V$ denote $rangeF_*$ and $(rangeF_*)^\perp$ part of $\nabla^N_U V,$ respectively.
The distribution $(rangeF_*)^\perp$ is totally geodesic if and only if  $R(\nabla^N_U V) = 0.$
Therefore, $(rangeF_*)^\perp$ is totally geodesic iff $$\nabla^N_U V = \nabla^{F^\perp}_U V.$$ Throughout this paper, we consider the map in which $(rangeF_*)^\perp$ is totally geodeic.
\begin{definition}\cite{Johnson}\label{deftgf}
Let $M$ be a Riemannian manifold and let $\mathcal{F}$ be a k- codimension foliation on $M.$ $\mathcal{F}$ is totally geodesic if each leaf $\mathcal{L}$ is a totally geodesic submanifold of $M,$ that is, any geodesic tangent to $\mathcal{L}$ at some point must lie within $\mathcal{L}.$
\end{definition}
\begin{theorem}({Frobenius Theorem}:)\cite{Lee book}\label{Frobenius}
A distribution is integrable if and only if it is involutive.         
\end{theorem}
\section{Clairaut semi-slant Riemnnian maps to K\"ahler manifolds}
In this section, we introduce the notion of a Clairaut semi-slant Riemannian map from a Riemannian manifold to a K\"ahler manifold. We investigate the geometry of this map.
\begin{definition}
Let $F:(M, g_M)\rightarrow (N, g_N, J)$ be a semi-slant Riemannian map from a Riemannian manifold to a K\"ahler manifold. Then, we say that $F$ is a Clairaut semi-slant Riemannian map if there exists a positive function $\tilde{s}$ on $N$ such that for any geodesic $\beta$ on $N,$ the function $(\tilde{s}\circ\beta)\sin\psi$ is constant, where $\psi(t)$ is the angle between $\dot{\beta}(t)$ and the horizontal subspace $\beta(t).$
\end{definition}
\begin{lemma}\label{geodesic}
Let $F:(M, g_M)\rightarrow (N, g_N, J)$ be a semi-slant Riemannian map from a Riemannian manifold
 to a K\"ahler manifold. If $\alpha$ is a geodesic on $(M, g_M),$ then the regular curve $\beta = F\circ \alpha$ is a geodesic on $N$ if and only if
\begin{equation}\label{hpart}
\begin{split}
&\cos^2\theta F_*(\nabla^M_X {X_2})+\cos^2\theta \nabla^N_{V}F_*{X_2}+F_*(\nabla^M_X {X_1})+\nabla ^N_{V}F_*{X_1}+S_{\omega(\phi F_*{X_2})}F_*X+\phi(S_{\omega F_*{X_2}}F_*X)\\&-B(\nabla^{F\perp}_{X}\omega F_*{X_2})-B(\nabla^{F\perp}_V \omega F_*{X_2})-B\big((\nabla F_*)(X, ^{*}F_*BV)\big)-\phi(F_*(\nabla^M_{X}{^{*}F_*BV}))\\&+\phi(S_{CV}F_*X)-B(\nabla^{F\perp}_X CV)-B(\nabla^{F\perp}_V CV)-\phi(\nabla^N_V BV) = 0 
\end{split}
\end{equation}
and
\begin{equation}\label{vpart}
\begin{split}
&\cos^2\theta(\nabla F_*)(X, X_2)+(\nabla F_*)(X, X_1)-\nabla^{F\perp}_{X}\omega(\phi F_*X_2)-\nabla^{F\perp}_V\omega(\phi F_*X_2)+\omega(S_{\omega F_*{X_2}}F_*X)\\&-C(\nabla^{F\perp}_X \omega F_*{X_2})-C(\nabla^{F\perp}_{V}\omega F_*{X_2})-C\big((\nabla F_*)(X, ^{*}F_*BV)\big)-\omega(F_*(\nabla^M_{X}{^{*}F_*BV}))\\&+\omega(S_{CV}F_*X)-C(\nabla^{F\perp}_X CV)-C(\nabla^{F\perp}_V CV)-\omega(\nabla^N_V BV)=0, 
\end{split}
\end{equation}
where $F_*X = F_*{X_1}+F_*{X_2},~~{X_1}\in\Gamma(D_1),$  ${X_2}\in\Gamma(D_2)$ and $V$ are the vertical and horizontal components of $\dot{\beta},$ respectively, $\nabla^N$ is the Levi–Civita connection on $N$ and $\nabla^{F\perp}$ is a linear connection on $(range F_*)^\perp.$
\end{lemma}
\begin{proof}
Let $\alpha$ be a geodesic on $M$ and $\beta = F\circ \alpha$ be a regular curve on $N.$ Suppose $F_*X$ and $V$ are the vertical and horizontal components respectively, of $\dot{\beta}$. Since $N$ is a K\"ahler manifold, we can write
\begin{equation*}
\nabla^N_{\dot{\beta}}\dot{\beta} = -J\nabla^N_{\dot{\beta}}J\dot{\beta}
\end{equation*} 
which gives
\begin{equation*}
\nabla^N_{\dot{\beta}}\dot{\beta} = -J\nabla^N_{F_*X+V}J(F_*X+V).
\end{equation*}
Since $F$ is a semi-slant Riemannian map, $F_*X = F_*{X_1}+F_*{X_2},$ where $F_*{X_1}\in\Gamma(F_*D_1)$ and $F_*{X_2}\in\Gamma(F_*D_2).$\\ Above equation can be rewritten as
\begin{equation*}
\begin{split}
\nabla^N_{\dot{\beta}}\dot{\beta} =-J\nabla^N_{F_*X+V}JF_*{X_1}-J\nabla^N_{F_*X+V}JF_*{X_2}-J\nabla^N_{F_*X+V}JV.  
\end{split}
\end{equation*}
Making use of \eqref{nablaJ}, \eqref{SJU} and \eqref{SJV} in the above equation, we obtain
\begin{equation}\label{nablaBeta}
\begin{split}
\nabla^N_{\dot{\beta}}\dot{\beta} = &\nabla^N_{F_*X}F_*{X_1}+\nabla^N_VF_*{X_1}-\nabla^N_{F_*X+V}(\phi^2F_*{X_2}+\omega(\phi F_*{X_2}))-J\nabla^N_{F_*X+V}\omega F_*X_2\\&-J\nabla^N_{F_*X+V}(BV+CV).
\end{split}
\end{equation}
With the help of \eqref{Sv} and \eqref{SFF}, we get
\begin{equation*}
\nabla^N_{F_*X}F_*{X_1} = (\nabla F_*)(X, X_1)+F_*(\nabla^M_X X_1),  
\end{equation*}
\begin{equation*}
\nabla^N_{F_*X}\omega F_*X_2 = -S_{\omega F_*{X_2}}F_*X+\nabla^{F\perp}_X\omega F_*X_2,
\end{equation*}
\begin{equation*}
\nabla^N_{F_*X}F_*{X_2} = (\nabla F_*)(X, X_2)+F_*(\nabla^M_X X_2),
\end{equation*}
\begin{equation*}
\nabla^N_{F_*X}BV = (\nabla F_*)(X, ^{*}F_*BV)+F_*(\nabla^M_{X}{^{*}F_*BV}),
\end{equation*}
\begin{equation*}
\nabla^N_{F_*X}CV = -S_{CV}F_*X+\nabla^{F\perp}_X CV,
\end{equation*}
\begin{equation*}
\nabla^N_{F_*X}{\omega(\phi F_*{X_2})} = -S_{\omega(\phi F_*{X_2})}F_*X+\nabla^{F\perp}_X\omega(\phi F_*{X_2}).
\end{equation*}
By the metric compatibility condition, we obtain
\begin{equation*}
\nabla^N_V F_*{X_1}\in\Gamma(rangeF_*),
\end{equation*}
\begin{equation*}
\nabla^N_V F_*{X_2}\in\Gamma(rangeF_*),
\end{equation*}
\begin{equation*}
\nabla^N_{V}BV\in\Gamma(range F_*).
\end{equation*}
Since $(rangeF_*)^\perp$ is totally geodesic, we get
\begin{equation*}
\nabla^N_V \omega(\phi F_*{X_2}) = \nabla^{F\perp}_V \omega(\phi F_*{X_2}),
\end{equation*}
\begin{equation*}
\nabla^N_V \omega F_*{X_2} = \nabla^{F\perp}_V \omega F_*{X_2},
\end{equation*}
\begin{equation*}
\nabla^N_V CV = \nabla^{F\perp}_V CV.
\end{equation*}
Making use of above equations in \eqref{nablaBeta}, we obtain
\begin{equation}
\begin{split}
\nabla^N_{\dot{\beta}}\dot{\beta} &=  (\nabla F_*)(X, X_1)+F_*(\nabla^M_X X_1)+\nabla^N_V F_*{X_1}+\cos^2\theta\big((\nabla F_*)(X, X_2)\\&+F_*(\nabla^M_{X}X_2)\big)+\cos^2\theta\nabla^N_V F_*X_2+S_{\omega(\phi F_*X_2)}F_*X-\nabla^{F^\perp}_X\omega(\phi F_*X_2)\\&-\nabla^{F^\perp}_V\omega(\phi F_*{X_2})+\phi(S_{\omega F_*{X_2}}F_*X)+\omega(S_{\omega F_*{X_2}}F_*X)-B(\nabla^{F^\perp}_X\omega F_*{X_2})\\&-C(\nabla^{F^\perp}_X\omega F_*{X_2})-B(\nabla^{F^\perp}_V \omega F_*{X_2})-C(\nabla^{F^\perp}_V \omega F_*{X_2})-B((\nabla F_*)(X, ^{*}F_*BV)\\&-C((\nabla F_*)(X, ^{*}F_*BV)-\phi(F_*(\nabla^M_{X}{^{*}F_*BV}))-\omega(F_*(\nabla^M_{X}{^{*}F_*BV}))\\&+\phi(S_{CV}F_*X)+\omega(S_{CV}F_*X)-B(\nabla^{F^\perp}_X CV)-C(\nabla^{F^\perp}_X CV)-B(\nabla^{F^\perp}_V CV)\\&-C(\nabla^{F^\perp}_V CV)-\phi(\nabla^N_V BV)-\omega(\nabla^N_V BV).
\end{split}
\end{equation}
We know that a regular curve $\beta$ is a geodesic on $N$ if and only if $\nabla^N_{\dot{\beta}} \dot{\beta} = 0$ this implies $\mathcal{V}\nabla^N_{\dot{\beta}} \dot{\beta} = 0$ and $\mathcal{H}\nabla^N_{\dot{\beta}} \dot{\beta} = 0.$ Taking the vertical and horizontal components of the vector field in the above equation, we get \eqref{hpart} and \eqref{vpart}, respectively.
\end{proof}
\begin{theorem}
Let $F$ be a semi-slant Riemannian map with connected fibers from a
Riemannian manifold $(M, g_M)$ to a K\"ahler manifold $(N, g_N, J)$ and consider $\beta$ to be geodesic on $N.$ Then $F$ is a Clairaut semi-slant Riemannian map with $\tilde{s} = e^g$ if and only if 
\begin{equation}\label{semi-slant NSC}
\begin{split}
\frac{d(g \circ \beta)}{dt} g_N(F_*X, F_*X)&= g_N\big(-B((\nabla F_*)(X, {^{*}F_* BV}))-\phi (F_*(\nabla^M_X {^{*}F_*BV})\\&+\phi(S_{CV}F_*X)-B(\nabla^{F^\perp}_X CV)-B(\nabla^{F^\perp}_V CV)\\&-\phi(\nabla^N_V CV), F_*X\big),
\end{split}
\end{equation}
where $F_*X, V$ are the vertical and horizontal components of $\dot{\beta},$ respectively.
\end{theorem}
\begin{proof}
Let $\beta$ be a geodesic on $N$ with constant speed $\sqrt{k}$, that is, $k = ||\dot{\beta}||^2.$
Let the vertical and horizontal parts of $\dot{\beta}$ be $F_*X$ and $V,$ respectively. Then, we get
\begin{equation}\label{gNFXFX}
g_N(F_*X, F_*X) = k\sin^2\psi(t), 
\end{equation}
\begin{equation*}\label{gNVV}
g_N(V, V) = k\cos^2\psi(t),
\end{equation*}
where $\psi(t)$ is the angle between $\dot{\beta}(t)$ and the horizontal subspace at $\beta(t).$\\
Differentiating \eqref{gNFXFX} with respect to `t', we get 
\begin{equation*}
\frac{d}{dt}g_N(F_*X, F_*X) = 2k \sin \psi(t) \cos \psi(t)\frac{d\psi}{dt}
\end{equation*}
which gives
\begin{equation}\label{semigNsincos}
g_N(\nabla^N_{\dot{\gamma}}F_*X, F_*X) = k \sin \psi(t) \cos \psi(t)\frac{d\psi}{dt}.
\end{equation}
Since $F$ is a semi-slant Riemannian map, $F_*X = F_*{X_1}+F_*{X_2},$ we have
\begin{equation*}\label{snsc}
g_N(\nabla^N_{\dot{\gamma}}F_*X, F_*X) = g_N(\nabla^N_{\dot{\gamma}}(F_*{X_1}+F_*{X_2}), F_*{X}).
\end{equation*}
From \eqref{semigNsincos} and the above equation, we have
\begin{equation*}
\begin{split}
k \sin \psi(t) \cos \psi(t)\frac{d\psi}{dt}& = g_N(\nabla^N_{F_*X}F_*{X_1}+\nabla^N_{F_*X}F_*{X_2}+\nabla^N_V F_*{X_1}\\&+\nabla^N_V F_*{X_2}, F_*X)
\end{split}
\end{equation*}
which gives
\begin{equation}
\begin{split}
k \sin \psi(t) \cos \psi(t)\frac{d\psi}{dt}& = g_N(\nabla^N_{F_*X}F_*{X_1}-J\nabla^N_{F_*X}JF_*{X_2}+\nabla^N_V F_*{X_1}\\&-J\nabla^N_V JF_*{X_2}, F_*X).
\end{split}
\end{equation}
With the help of \eqref{Sv}, \eqref{SFF} and the above equation, we get
\begin{equation*}
\begin{split}
k \sin \psi(t) \cos \psi(t)\frac{d\psi}{dt} &= g_N(F_*(\nabla^M_X {X_1})+\nabla^N_V F_*{X_1}+\cos^2\theta F_*(\nabla^M_X{X_2})+S_{\omega(\phi F_*{X_2})}F_*X\\&+\phi(S_{\omega F_*{X_2}}F_*X)-B(\nabla^{F^\perp}_X \omega F_*{X_2})+\cos^2\theta\nabla^N_V F_*{X_2}\\&-B(\nabla^{F^\perp}_V \omega F_*{X_2}), F_*X).
\end{split}
\end{equation*}
Making use of Lemma \ref{geodesic} in the above equation, we get
\begin{equation}\label{sincos}
\begin{split}
k \sin \psi(t) \cos \psi(t)\frac{d\psi}{dt}& = g_N\big(B((\nabla F_*)(X, {^{*}F_* BV}))+\phi (F_*(\nabla^M_X {^{*}F_*BV})\\&-\phi(S_{CV}F_*X)+B(\nabla^{F^\perp}_X CV)+B(\nabla^{F^\perp}_V CV)\\&+\phi(\nabla^{F^\perp}_V CV), F_*X\big).
\end{split}
\end{equation}
Further, $F$ is a Clairaut semi-slant Riemannian map with $\tilde{s} = e^g$ if and only if  
\begin{equation*}
\frac{d}{dt}( e^{g \circ \beta} \sin \psi(t)) = 0
\end{equation*}
which implies  
\begin{equation*}
e^{g\circ\beta}\cos\psi(t)\frac{d\psi}{dt}+e^{g\circ\beta}\frac{d(g\circ \beta)}{dt}\sin\psi(t) = 0.
\end{equation*}
Since $e^{g \circ \beta}$ is a positive function, we get 
\begin{equation*}
\frac{d(g \circ \beta)}{dt} \sin\psi(t) + \cos\psi(t) \frac{d\psi}{dt} = 0.
\end{equation*}
Multiplying the above equation with $k\sin\psi$ and comparing with \eqref{sincos}, we get \eqref{semi-slant NSC}.
\end{proof}
\begin{lemma}\cite{MeenaCRMap}\label{harmonic}
Let $F:(M, g_M)\rightarrow (N, g_N)$ be a Clairaut Riemannian map with $\tilde{s} = e^g$ between Riemannian manifolds such that $kerF_*$ is minimal. Then $F$ is harmonic if and only if $g$ is constant.
\end{lemma}
\begin{theorem}
Let $F:(M, g_M)\rightarrow (N, g_N, J)$ be a Clairaut semi-slant Riemannian map with $\tilde{s} = e^g$ from a Riemannian manifold to a K\"ahler manifold such that \text{ker}$F_*$ is minimal. Then for $F_*X = F_*{X_1}+F_*{X_2}$ with $F_*{X_1}\in\Gamma(F_*{D_1})$ and $F_*{X_2}\in\Gamma(F_*{D_2}),$ $F$ is a harmonic map if and only if
\begin{equation}
trace\Big(\nabla^{F^\perp}_X\omega(\phi F_*{X_2})+C(\nabla^{F^\perp}_X \omega F_*{X_2})\Big) = 0.
\end{equation}
\end{theorem}
\begin{proof}
We know that
\begin{equation*}
(\nabla F_*)(X, X) = \nabla^N_{F_{*}X}F_*X - F_{*}(\nabla^M_X X).
\end{equation*}
We can write $F_*X = F_*{X_1}+F_*{X_2},$ where $F_*{X_1}\in\Gamma(F_*{D_1})$ and $F_*{X_2}\in\Gamma(F_*{D_2})$ and making use of \eqref{nablaJ} in the above equation, we get
\begin{equation*}
\begin{split}
(\nabla F_*)(X, X) = \nabla^N_{F_{*}X}F_*{X_1}-J\nabla^N_{F_{*}X}JF_*{X_2}- F_{*}(\nabla^M_X X).
\end{split}
\end{equation*}
Using  \eqref{Sv}, \eqref{SFF}, \eqref{SJU} and \eqref{SJV} in the above equation, we obtain
\begin{equation*}
\begin{split}
(\nabla F_*)(X, X)& = (\nabla F_*)(X, X_1)+F_*(\nabla^M_X {X_1})+\cos^2\theta(\nabla F_*)(X, X_2)+\cos^2\theta F_*(\nabla^M_X {X_2})\\&+S_{\omega(\phi F_*{X_2})}F_*X-\nabla^{F^\perp}_X \omega(\phi F_*{X_2})+\phi(S_{\omega F_*{X_2}}F_*X)+\omega(S_{\omega F_*{X_2}}F_*X)\\&-B(\nabla^{F^\perp}_X \omega F_*{X_2})-C(\nabla^{F^\perp}_X \omega F_*{X_2})-F_*(\nabla^M_X X).
\end{split}
\end{equation*}
Since the range part of second fundamental form is zero, equating the $(range F_*)^\perp$ part of above equation and using Theorem \ref{NSC}, we get
\begin{equation*}
\begin{split}
-g_M(X, X)\nabla^N g =& -g_M(X, X_1)\nabla^N g-\cos^2\theta g_M(X, X_2)\nabla^N g-\nabla^{F^\perp}_X\omega(\phi F_*{X_2})\\&-\omega((\omega F_*{X_2})(g)F_*X)-C(\nabla^{F\perp}_X \omega F_*{X_2}).
\end{split}
\end{equation*}
Making use of Lemma \ref{harmonic} and taking trace in the above equation, we get the required result.
\end{proof}
 \begin{theorem}
 Let $F:(M, g_M)\rightarrow (N, g_N, J)$ be a Clairaut semi-slant Riemannian map with $\tilde{s}=e^g$ from a Riemannian manifold to a K\"ahler manifold. Then, for any $F_*X, F_*Y\in\Gamma(F_*{D_1}),$ $F_*Z\in\Gamma(F_*{D_2})$ and $V\in\Gamma(rangeF_*)^\perp,$ the following assertions are equivalent:
 \begin{itemize}
 \item [(i)] The invariant distribution $F_*(D_1)$ is a totally geodesic foliation on $N,$
 \item [(ii)]$g_N\big(F_*(\nabla^M_X {^{*}F_*(JF_*Y)}), \phi F_*Z\big)-g_M(X, {^{*}F_*(JF_*Y)})g_N\big(\nabla^N g , \omega F_*Z\big) = 0$\\ and\\ $g_N\big(F_*(\nabla^M_X {^{*}F_*JF_*Y}), BV\big)-CV(g)g_N(F_*X, JF_*Y) = 0,$ 
 \item [(iii)]
$g_M(X, Y)g_N\big(B\nabla^N g+C\nabla^N g, \phi F_*Z+\omega F_*Z) = \cos^2\theta g_N(F_*(\nabla^M_X Y), F_*Z)\\-g_N(F_*(\nabla^M_X Y), B(\omega F_*Z)).$
\end{itemize}
\end{theorem}
 \begin{proof}
 For any $F_*X, F_*Y\in\Gamma(F_*{D_1})$ and $F_*Z\in\Gamma(F_*{D_2}),$ we have
 \begin{equation*}
 \begin{split}
 g_N(\nabla^N_{F_*X}F_*Y, F_*Z) = g_N(\nabla^N_{F_*X}JF_*Y, JF_*Z).
 \end{split}
 \end{equation*}
 Making use of \eqref{SFF} and \eqref{SJU} in the above equation, we obtain
 \begin{equation*}
 \begin{split}
 g_N(\nabla^N_{F_*X}F_*Y, F_*Z) & = g_N\big(F_*(\nabla^M_X {^{*}F_*(JF_*Y)}), \phi F_*Z\big)\\&+g_N\big((\nabla F_*)(X, {^{*}F_*(JF_*Y)}) , \omega F_*Z\big).
 \end{split}
 \end{equation*}
 Using Theorem \ref{NSC} in the above equation, we get
 \begin{equation}\label{sdtgf}
 \begin{split}
  g_N(\nabla^N_{F_*X}F_*Y, F_*Z) & = g_N\big(F_*(\nabla^M_X {^{*}F_*(JF_*Y)}), \phi F_*Z\big)\\&-g_M(X, {^{*}F_*(JF_*Y)})g_N\big(\nabla^N g , \omega F_*Z\big).
 \end{split}
 \end{equation}
 On the other hand, for $V\in\Gamma(range F_*)^\perp,$ we have
 \begin{equation*}
 g_N(\nabla^N_{F_*X}F_*Y, V) = g_N(\nabla^N_{F_*X}JF_*Y, JV).
 \end{equation*}
 Making use of \eqref{Sv} and \eqref{SFF} in the above equation, we get
 \begin{equation*}
g_N(\nabla^N_{F_*X}F_*Y, V) = g_N\big(F_*(\nabla^M_X {^{*}F_*JF_*Y}), BV\big)+g_N(S_{CV}F_*X, JF_*Y).
 \end{equation*}
 Using Theorem \ref{NSC} in the above equation, we obtain
\begin{equation*}
 g_N(\nabla^N_{F_*X}F_*Y, V) = g_N\big(F_*(\nabla^M_X {^{*}F_*JF_*Y}), BV\big)-CV(g)g_N(F_*X, JF_*Y). \end{equation*}
 From \eqref{sdtgf}, above equation and using Definition \ref{deftgf}, the proof of $(i)\iff (ii)$ follows.\\
 For any $F_*X, F_*Y\in\Gamma(F_*{D_1})$ and $F_*Z\in\Gamma(F_*{D_2}),$ we have 
 \begin{equation*}
 g_N(\nabla^N_{F_*X}F_*Y, F_*Z) = g_N(J\nabla^N_{F_*X}F_*Y, JF_*Z).
 \end{equation*}
 Using \eqref{SFF}, \eqref{SJU} and \eqref{SJV} in the above equation, we obtain
\begin{equation*}
\begin{split}
g_N(\nabla^N_{F_*X}F_*Y, F_*Z)& =  g_N\big(B(\nabla F_*)(X, Y)+C(\nabla F_*)(X, Y), \phi F_*Z+\omega F_*Z\big)\\&-g_N\big(F_*(\nabla^M_X Y), J(\phi F_*Z+\omega F_*Z)\big). 
\end{split}
\end{equation*}
Making use of Theorem \ref{NSC} in the above equation, we get
\begin{equation*}
\begin{split}
g_N(\nabla^N_{F_*X}F_*Y, F_*Z)& = -g_M(X, Y)g_N\big(B\nabla^N g+C\nabla^N g, \phi F_*Z+\omega F_*Z)\\&+\cos^2\theta g_N(F_*(\nabla^M_X Y), F_*Z)-g_N(F_*(\nabla^M_X Y), B\omega F_*Z).
\end{split}
\end{equation*}
From above equation and using Definition \ref{deftgf}, the proof of $(i)\iff (iii)$ follows.
\end{proof}
\begin{theorem}
Let $F:(M, g_M)\rightarrow (N, g_N, J)$ be a Clairaut semi-slant Riemannian map with $\tilde{s}=e^g$ from a Riemannian manifold to a K\"ahler manifold. Then, for any $F_*X, F_*Y\in\Gamma(F_*{D_2}),$ $F_*Z\in\Gamma(F_*{D_1})$ and $V\in\Gamma(rangeF_*)^\perp,$ the following assertions are equivalent:
\begin{itemize}
\item [(i)] The slant distribution $F_*(D_2)$ is a totally geodesic foliation on $N,$
\item [(ii)]$F_*(\nabla^M_X {^{*}F_*\phi F_*Y})\in\Gamma(F_*D_2)$\\and\\ $g_M(X, {^{*}F_*\phi F_*Y})g_N(C(\nabla^N g), V)-g_N(\omega F_*(\nabla^M_X {^{*}F_*\phi F_*Y}), V)-\omega F_*Y(g)g_N(\omega F_*X, V)-C(\nabla^{F\perp}_X \omega F_*Y), V) = 0.$
\end{itemize}
\end{theorem}
\begin{proof}
For any $F_*X, F_*Y\in\Gamma(F_*{D_2})$ and $F_*Z\in\Gamma(F_*{D_1}),$ we have
\begin{equation*}
\begin{split}
g_N(\nabla^N_{F_*X}F_*Y, F_*Z) = g_N(\nabla^N_{F_*X}JF_*Y, JF_*Z).
\end{split}
\end{equation*}
Making use of \eqref{SFF} and \eqref{SJU} in the above equation, we get
\begin{equation*}
g_N(\nabla^N_{F_*X}F_*Y, F_*Z) = g_N(F_*(\nabla^M_X {^{*}F_*\phi F_*Y}), JF_*Z)-g_N(S_{\omega F_*Y}F_*X, JF_*Z).
\end{equation*}
Using Theorem \ref{NSC} in the above equation, we obtain
\begin{equation}\label{sstgf}
 g_N(\nabla^N_{F_*X}F_*Y, F_*Z)  = g_N(F_*(\nabla^M_X {^{*}F_*\phi F_*Y}), JF_*Z)   
\end{equation}
On the other hand, for $V\in\Gamma(range F_*)^\perp,$ we have
\begin{equation*}
g_N(\nabla^N_{F_*X}F_*Y, V) =  g_N(\nabla^N_{F_*X}JF_*Y, JV).
\end{equation*}
Using \eqref{Sv} and \eqref{SFF} in the above equation, we get
\begin{equation*}
\begin{split}
g_N(\nabla^N_{F_*X}F_*Y, V) &= -g_N\big(J(\nabla F_*)(X, {^{*}F_*\phi F_*Y})+JF_*(\nabla^M_X {^{*}F_*\phi F_*Y})\\&+J(-S_{\omega F_*Y}F_*X+\nabla^{F\perp}_X \omega F_*Y), V\big).
\end{split}
\end{equation*}
Making use of \eqref{SJU}, \eqref{SJV} and Theorem \ref{NSC} in the above equation, we get
\begin{equation*}
\begin{split}
g_N(\nabla^N_{F_*X}F_*Y, V) &= g_M(X, {^{*}F_*\phi F_*Y})g_N(C(\nabla^N g), V)-g_N(\omega F_*(\nabla^M_X {^{*}F_*\phi F_*Y}), V)\\&-\omega F_*Y(g)g_N(\omega F_*X, V)-g_N(C(\nabla^{F\perp}_X \omega F_*Y), V)
\end{split}
\end{equation*}
From the above equation, \eqref{sstgf} and Definition \ref{deftgf}, the proof of $(i)\iff (ii)$ follows.
\end{proof}
From above two theorems, we have the following decomposition Theorem.
\begin{theorem}
Let $F:(M, g_M)\rightarrow (N, g_N, J)$ be a Clairaut semi-slant Riemannian map with $\tilde{s}=e^g$ from a Riemannian manifold to a K\"ahler manifold with integrable distribution $(range F_*).$ Then for $F_*X, F_*Y \in\Gamma(range F_*), F_*Z\in\Gamma(F_*D_2)$ and $V\in\Gamma(range F_*)^\perp,$ the leaf of $(range F_*)$ is a locally product Riemannian manifold of $F_*(D_1)$ and $F_*(D_2)$ if and only if
\begin{equation*}
g_M(X, Y)g_N\big(B\nabla^N g+C\nabla^N g, \phi F_*Z+\omega F_*Z)= \cos^2\theta g_N(F_*(\nabla^M_X Y), F_*Z)-g_N(F_*(\nabla^M_X Y), B\omega F_*Z),
\end{equation*}
$F_*(\nabla^M_X {^{*}F_*\phi F_*Y})\in\Gamma(F_*D_2),$ and
\begin{equation*}
\begin{split}
& g_M(X, {^{*}F_*\phi F_*Y})g_N(C(\nabla^N g), V)-g_N(\omega F_*(\nabla^M_X {^{*}F_*\phi F_*Y}), V)-\omega F_*Y(g)g_N(\omega F_*X, V)\\&-C(\nabla^{F\perp}_X \omega F_*Y), V) = 0.  
 \end{split}
\end{equation*}
\end{theorem}
\begin{theorem}
Let $F:(M, g_M)\rightarrow (N, g_N, J)$ be a Clairaut semi-slant Riemannian map with $\tilde{s} = e^g$ from a Riemannian manifold to a K\"ahler manifold. Then, for $F_*X, F_*Y\in\Gamma(F_*D_2), F_*Z\in\Gamma (F_*D_1),$ the distribution $F_*(D_2)$ is a integrable if and only if
\begin{equation}
\begin{split}
&g_M(X, {^{*}F_*\phi F_*Y})g_N(B\nabla^N g, F_*Z)-g_N(\phi F_*(\nabla^M_X {^{*}F_*\phi F_*Y}), F_*Z)\\&+g_M(F_*(\nabla^M_Y {^{*}F_*\phi F_*X}), JF_*Z) = 0.
\end{split}
\end{equation}
\end{theorem}
\begin{proof}
For $F_*X, F_*Y\in\Gamma (F_*D_2)$ and $F_*Z\in\Gamma (F_*D_1),$ we have
\begin{equation*}
\begin{split}
g_N([F_*X, F_*Y], F_*Z) &= g_N(\nabla^N_{F_*X}F_*Y-\nabla^N_{F_*Y}F_*X, F_*Z)\\&=g_N(\nabla^N_{F_*X}JF_*Y, JF_*Z)-g_N(\nabla^N_{F_*Y}JF_*X, JF_*Z).
\end{split}
\end{equation*}
Using \eqref{Sv}, \eqref{SFF} and \eqref{SJU} in the above equation, we get
\begin{equation*}
\begin{split}
g_N([F_*X, F_*Y], F_*Z) &= -g_N\Big(J\big((\nabla F_*)(X, {^{*}F_*\phi F_*Y})+F_*(\nabla^M_X{^{*}F_*\phi F_*Y})\big), F_*Z\Big)\\&-g_N(S_{\omega F_*Y}F_*X, JF_*Z)-g_N\Big((\nabla F_*)(Y, {^{*}F_*\phi F_*X})\\&+F_*(\nabla^M_Y{^{*}F_*\phi F_*X})-S_{\omega F_*X}F_*Y, JF_*Z\Big)
\end{split}
\end{equation*}
Making use of \eqref{SJU}, \eqref{SJV} and Theorem \ref{NSC} in the above equation, we get
\begin{equation*}
\begin{split}
g_N([F_*X, F_*Y], F_*Z) &= g_M(X, {^{*}F_*\phi F_*Y})g_N(B\nabla^N g, F_*Z)-g_N(\phi F_*(\nabla^M_X {^{*}F_*\phi F_*Y}), F_*Z)\\&+g_M(F_*(\nabla^M_Y {^{*}F_*\phi F_*X}), JF_*Z)
\end{split}
\end{equation*}
We get the required result from Theorem \ref{Frobenius} and the above equation.
\end{proof}
\begin{theorem}
Let $F : (M, g_M) \rightarrow (N, g_N, J)$ be a Clairaut semi-slant Riemannian map with $\tilde{s} = e^g$ from a Riemannian manifold to a K\"ahler manifold. Then $F$ is totally geodesic if and only if the following conditions are satisfied:
\begin{itemize}
\item [(i)]$kerF_*$ and  $(kerF_*)^\perp$ are totally geodesic distributions,
\item [(ii)] For $X\in\Gamma(kerF_*)^\perp$ and $Y\in\Gamma(D_1)$ with $JF_*Y = F_*Z,$ we have 
\begin{equation*}
\begin{split}
g_M(X, Z)\big(B(\nabla^N g)+C(\nabla ^N g)\big)-\phi F_*(\nabla^M_X Z)-\omega F_*(\nabla^M_X Z)-F_*(\nabla^M_X Y) = 0,  
\end{split}
\end{equation*}
\item [(iii)] For $X, Y\in\Gamma(D_2),$ we have 
\begin{equation*}
\begin{split}
&\cos^2\theta \nabla^N_{F_*X}F_*Y+\omega(\phi F_*Y)(g)F_*X+\nabla^{F\perp}_X \omega(\phi F_*Y)+\omega F_*Y(g)(\phi F_*X + \omega F_*X)\\&+B(\nabla^{F\perp}_X \omega F_*Y)+C(\nabla^{F\perp}_X \omega F_*Y)-F_*(\nabla^M_X Y) = 0.   
\end{split}    
\end{equation*} 
\end{itemize}
\end{theorem}
\begin{proof}
We know that a Riemannian map $F$ is totally geodesic if and only if $(\nabla F_*) = 0,$ that is, 
\begin{equation}
(\nabla F_*)(U, V) = 0,   
\end{equation}
\begin{equation}
(\nabla F_*)(U, X) = 0, 
\end{equation}
and
\begin{equation}\label{semitgm}
(\nabla F_*)(X, Y) = 0.
\end{equation}
$\forall U, V\in\Gamma (kerF_*)$ and $\forall X, Y\in\Gamma(kerF_*)^\perp.$
The first two equations implies that $kerF_*$ and $(kerF_*)^\perp$ are totally geodesic which proves (i) and (ii). For the proof of (iii) and (iv), we proceed as follows:\\
For $X\in\Gamma(kerF_*)^\perp$ and $Y\in\Gamma(D_1),$ we have
\begin{equation*}
(\nabla F_*)(X, Y) = \nabla^N_{F_*X}F_*Y-F_*(\nabla^M_X Y).
\end{equation*}
\begin{equation*}
\begin{split}
(\nabla F_*)(X, Y) &= J\nabla^N_{F_*X}JF_*Y-F_*(\nabla^M_X Y),\\&= J\nabla^N_{F_*X}F_*Z- F_*(\nabla^M_X Y),
\end{split}
\end{equation*}
where $JF_*Y = F_*Z$ (say). 
Using  \eqref{SFF} and Theorem \ref{NSC} in the above equation, we obtain
\begin{equation*}
\begin{split}
(\nabla F_*)(X, Y) = J\big(-g_M(X, Z)\nabla^N g+F_*(\nabla^M_X Z)\big)-F_*(\nabla^M_X Y).
\end{split}
\end{equation*}
Making use of \eqref{SJU} and \eqref{SJV} in the above equation, we get
\begin{equation*}
\begin{split}
 (\nabla F_*)(X, Y) &= -g_M(X, Z)\big(B(\nabla^N g)+C(\nabla ^N g)\big)+\phi F_*(\nabla^M_X Z)+\omega F_*(\nabla^M_X Z)\\&-F_*(\nabla^M_X Y).  
 \end{split}
\end{equation*}
From \eqref{semitgm} and the above equation, we get the proof of part (ii).\\
Now, for $X, Y\in\Gamma(D_2),$ we have
\begin{equation*}
\begin{split}
(\nabla F_*)(X, Y) &= -J\nabla^N_{F_*X}JF_*Y-F_*(\nabla^M_X Y)\\& = -J\nabla^N_{F_*X}(\phi F_*Y+\omega F_*Y)-F_*(\nabla^M_X Y)\\& = -\nabla^N_{F_*X}\big(\phi^2 F_*Y+\omega(\phi F_*Y)\big)-J(-S_{\omega F_*Y}F_*X+\nabla^{F\perp}_X \omega F_*Y)-F_*(\nabla^M_X Y)\\& = \cos^2\theta\nabla^N_{F_*X}F_*Y+S_{\omega(\phi F_*Y)}F_*X-\nabla^{F\perp}_X \omega(\phi F_*Y)+\phi(S_{\omega F_*Y}F_*X)\\&+ \omega (S_{\omega F_*Y}F_*X)-B(\nabla^{F\perp}_X \omega F_*Y)-C(\nabla^{F^\perp}_X \omega F_*Y)-F_*(\nabla^M_X Y).
\end{split}
\end{equation*}
Making use of Theorem \ref{NSC} in the above equation, we obtain
\begin{equation*}
\begin{split}
(\nabla F_*)(X, Y) &= \cos^2\theta\nabla^N_{F_*X}F_*Y-\omega(\phi F_*Y)(g)F_*X-\nabla^{F^\perp}_X \omega(\phi F_*Y)-C(\nabla^{F^\perp}_X \omega F_*Y)\\&-\omega F_*Y(g)(\phi F_*X + \omega F_*X)-B(\nabla^{F^\perp}_X \omega F_*Y)-F_*(\nabla^M_X Y).
\end{split}
\end{equation*} 
Making use of \eqref{semitgm} in the above equation, we get the proof of part (iii). 
\begin{example}
Let $M = N = \mathbb{R}^6$ be Euclidean spaces with Riemannian metrics $g_M = dx_1^2+dx_2^2+e^{2x_3}dx_3^2+e^{2x_3}dx_4^2+dx_5^2+dx_6^2$ and $g_N =  dy_1^2+dy_2^2+e^{2x_3}dy_3^2+e^{2x_3}dy_4^2+dy_5^2+dy_6^2.$
Define a complex structure $J$ on $N$ as $$J(y_1, y_2, y_3, y_4, y_5, y_6) = (-y_2, y_1, -y_4, y_3, -y_6, y_5)$$ and a smooth map $F:(M, g_M)\rightarrow (N, g_N, J)$ as
$$F(x_1, x_2, x_3, x_4, x_5, x_6) = (x_2\cos \alpha, 0, x_3, x_4, x_5, x_2\sin\alpha).$$\\Then, we have
$$\ker F_* = \text{span}\{e_1, e_6\},$$
$$(kerF_*)^\perp = \text{span}\{X_1= e_2, X_2= e_3, X_3= e_4, X_4= e_5\},$$
where $D_1 = \{e_3, e_4\}$ and $D_2 = \{e_2, e_5\},$
$$range F_* = \text{span}\{F_*{X_1} = \cos\alpha e_1'+\sin\alpha e_6', F_*{X_2} = e_3', F_*{X_3} = e_4', F_*{X_4} = e_5'\},$$
$$(rangeF_*)^\perp = \text{span}\{V_1= e_2', V_2 = -\sin\alpha e_1'+\cos\alpha e_6'\},$$
where $\{e_1 = \frac{\partial}{\partial x_1}, e_2 = \frac{\partial}{\partial x_2}, e_3 = e^{-x_3}\frac{\partial}{\partial x_3}, e_4 = e^{-x_3}\frac{\partial}{\partial x_4}, e_5 = \frac{\partial}{\partial x_5}, e_6 = \frac{\partial}{\partial x_6}\}$ and $\{e_1' = \frac{\partial}{\partial y_1}, e_2' = \frac{\partial}{\partial y_2}, e_3' = e^{-x_3}\frac{\partial}{\partial y_3}, e_4' = e^{-x_3}\frac{\partial}{\partial y_4}, e_5' = \frac{\partial}{\partial y_5}, e_6' = \frac{\partial}{\partial y_6}\}$ are bases of $T_pM$  and $T_{F(p)}N,$ respectively.\\ The condition for a smooth map between Riemannian manifolds to be a Riemannian map, that is, $g_M(X, X) = g_N(F_*X, F_*X)$ holds here. Therefore, $F$ is a Riemannian map from a Riemannian manifold to an almost Hermitian manifold.\\ Also $J(F_*X_1) = \cos\alpha e_2'-\sin\alpha e_5', J(F_*X_2) = e_4', J(F_*X_3)= -e_3', J(F_*X_4) = e_6', J(V_1) = -e_1', J(V_2) = -\sin\alpha e_2'-\cos\alpha e_5'.$ Clearly $F$ is a semi-slant Riemannian map with semi-slant angle $\alpha.$ The non-zero Christoffel symbols are $\Gamma_{3 3}^{3} = 1, \Gamma_{4 4}^3 = -1, \Gamma_{3 4}^4 = 1 = \Gamma_{4 3}^4.$
For any $\Tilde{X} = a_1e_1'+ a_2e_2'+ a_3e_3'+ a_4e_4'+ a_5e_5'+ a_6e_6'$ and $\Tilde{Y} = b_1e_1'+ b_2e_2'+ b_3e_3'+ b_4e_4'+ b_5e_5'+ b_6e_6',$ where $a_i, b_i\in\mathbb{R},$ we have $\nabla^N_{\Tilde{X}}J\Tilde{Y} = J\nabla^N_{\Tilde{X}}\Tilde{Y}$ which shows that the K\"ahlerian condition holds. Further, For $X\in(kerF_*)^\perp,$
$$g_M(X, X) = g_M(a_1e_2+a_2e_3+a_3e_4+a_4e_5, a_1e_2+a_2e_3+a_3e_4+a_4e_5)= a_1^2+a_2^2+a_3^2+a_4^2.$$
and $$(\nabla F_*)(X, X) = \mu_1e_2'+\mu_2(-\sin\alpha e_1'+\cos\alpha e_6')~~ for ~~ \mu_i\in\mathbb{R}.$$ Therefore, $(\nabla F_*)(X, X) = -g_M(X, X)\nabla^N g$ for $$g = \frac{\mu_1y_2+\mu_2(-\sin\alpha y_1+\cos\alpha y_6)}{a_1^2+a_2^2+a_3^2+a_4^2},$$ where atleast one $a_i$ is non-zero. This shows that $F$ is a Clairaut semi-slant Riemannian map.
\end{example}
\end{proof}
\section{Clairaut hemi-slant Riemannian maps to K\"ahler manifolds}
In this section, we define the Clairaut hemi-slant Riemannian map to K\"ahler manifold and investigate the geometry of this map.
\begin{definition}
Let $F:(M, g_M)\rightarrow (N, g_N, J)$ be a hemi-slant Riemannian map from a Riemannian manifold to a K\"ahler manifold. Then, we say that $F$ is a Clairaut hemi-slant Riemannian map if there exist a positive function $\tilde{s}$ on $N$ such that for any geodesic $\beta$ on $N,$ the function $(\tilde{s}\circ\beta)\sin\psi$ is constant, where $\psi(t)$ is the angle between $\dot{\beta}(t)$ and the horizontal subspace at $\beta(t).$
\end{definition}
\begin{theorem}\label{hemi-geodesic}
Let $F:(M, g_M)\rightarrow (N, g_N, J)$ be a hemi-slant Riemannian map from a Riemannian manifold to a K\"ahler manifold. If $\alpha$ is a geodesic on $(M, g_M),$ then the regular curve $\beta = F\circ \alpha$ is a geodesic on $N$ if and only if
\begin{equation}\label{hemiv}
\begin{split}
 &\phi(S_W F_*X)-B(\nabla^{F\perp}_X W)-B(\nabla^{F\perp}_V W)+\cos^2\theta F_*(\nabla^M_X {X_2})+\cos^2\theta\nabla^N_VF_*{X_2}\\&+S_{\omega(\phi F_*{X_2})}F_*X+\phi(S_{\omega F_*{X_2}}F_*X)-B(\nabla^{F\perp}_X \omega F_*{X_2})-B(\nabla^{F\perp}_V \omega F_*{X_2})\\&-B(\nabla F_*)(X, ^{*}F_*BV)-\phi(F_*(\nabla^M_X{^{*}F_*BV})+\phi(S_{CV}F_*X)-B(\nabla^{F\perp}_X CV)\\&-B(\nabla^{F\perp}_V CV)-\phi(\nabla^N_V BV) = 0,
\end{split}
\end{equation}
\begin{equation}\label{hemih}
\begin{split}
&\omega(S_W F_*X)-C(\nabla^{F\perp}_X W)-C(\nabla^{F\perp}_V W)+\cos^2\theta(\nabla F_*)(X, X_2)-\nabla^{F^\perp}_X \omega(\phi F_*{X_2})\\&-\nabla^{F\perp}_V\omega(\phi F_*{X_2})+\omega(S_{\omega F_*{X_2}}F_*X)-C(\nabla^{F\perp}_X \omega F_*{X_2})- C(\nabla^{F\perp}_V \omega F_*{X_2})\\&-C\big((\nabla F_*)(X, {^{*}F_*BV})\big)-\omega (F_*(\nabla^M_X {^{*}F_*BV}))+\omega (S_{CV}F_*X)-C(\nabla^{F\perp}_X CV)\\&-C(\nabla^{F\perp}_V CV)-\omega(\nabla^N_V BV) = 0,
\end{split}
\end{equation}
where $F_*X = F_*{X_1}+F_*{X_2}, JF_*{X_1} = W,~~F_*{X_1}\in\Gamma(D^\perp), F_*{X_2}\in\Gamma(D^\psi)$ and $V$ are the vertical and horizontal components of $\dot{\beta},$ respectively.
\end{theorem}
\begin{proof}
Let $\alpha$ be a geodesic on $M$ and $\beta = F\circ \alpha$ be a regular curve on $N.$ Consider $F_*X$ and $V$ are the vertical and horizontal components of $\dot{\beta},$ respectively. By the property of K\"ahler manifold, we can write
\begin{equation*}
\nabla^N_{\dot{\beta}}\dot{\beta} = -J\nabla^N_{\dot{\beta}}J\dot{\beta}.
\end{equation*}
Since $F$ is a hemi-slant Riemannian map, $F_*X = F_*{X_1}+F_*{X_2},$ where $F_*{X_1}\in\Gamma(D^\perp)$ and $F_*{X_2}\in\Gamma(D^{\psi}).$ Therefore, the above equation turns into
\begin{equation*}
\nabla^N_{\dot{\beta}}\dot{\beta} = -J\nabla^N_{\dot{\beta}}J(F_*{X_1}+F_*{X_2}+V)
\end{equation*}
which implies
\begin{equation*}
\nabla^N_{\dot{\beta}}\dot{\beta} = -J\nabla^N_{F_*X+V}W-J\nabla^N_{F_*X+V}JF_*{X_2}-J\nabla^N_{F_*X+V}JV,
\end{equation*}
where $JF_*{X_1} = W$ (say). Using \eqref{Sv}, \eqref{SFF}, \eqref{hemiJV} and \eqref{hemiJxi} in the above equation, we obtain
\begin{equation*}
\begin{split}
\nabla^N_{\dot{\beta}}\dot{\beta}& =\phi(S_W F_*X)+\omega(S_W F_*X)-B(\nabla^{F^\perp}_X W)-C(\nabla^{F^\perp}_X W)-B(\nabla^{F^\perp}_V W)-C(\nabla^{F^\perp}_V W)\\&+\cos^2\theta\big((\nabla F_*)(X, X_2)+F_*(\nabla^M_X X_2)\big)+\cos^2\theta \nabla^N_V F_*{X_2}+S_{\omega(\phi F_*{X_2})}F_*X-\nabla^{F^\perp}_X \omega(\phi F_*{X_2})\\&-\nabla^{F^\perp}_V \omega(\phi F_*{X_2})+\phi(S_{\omega F_*{X_2}}F_*X)+\omega(S_{\omega F_*{X_2}}F_*X)-B(\nabla^{F^\perp}_X \omega F_*{X_2})-C(\nabla^{F^\perp}_X \omega F_*{X_2})\\&-B(\nabla^{F^\perp}_V \omega F_*{X_2})- C(\nabla^{F^\perp}_V \omega F_*{X_2})-B\big((\nabla F_*)(X, ^{*}F_*BV)\big)- C\big((\nabla F_*)(X, ^{*}F_*BV)\big)\\&-\phi(F_*(\nabla^M_{X}{^{*}F_*BV})\big)-\omega(F_*(\nabla^M_{X}{^{*}F_*BV}))+\phi(S_{CV}F_*X)+\omega(S_{CV}F_*X)-B(\nabla^{F^\perp}_X CV)\\&-C(\nabla^{F^\perp}_X CV)-B(\nabla^{F^\perp}_V CV)-C(\nabla^{F^\perp}_V CV)-\phi(\nabla^N_V BV)-\omega(\nabla^N_V BV).
\end{split}
\end{equation*}
We know that a regular curve $\beta$ is a geodesic on $N$ if and only if $\nabla^N_{\dot{\beta}} \dot{\beta} = 0$ which implies $\mathcal{V}\nabla^N_{\dot{\beta}} \dot{\beta} = 0$ and $\mathcal{H}\nabla^N_{\dot{\beta}} \dot{\beta} = 0.$ Taking the vertical and horizontal components of the vector field in the above equation, we get \eqref{hemiv} and \eqref{hemih}.
\end{proof}
\begin{theorem}
Let $F:(M, g_M)\rightarrow (N, g_N, J)$ be a hemi-slant Riemannian map from a Riemannian manifold to a K\"ahler manifold and $\beta$ is a geodesic on $N$. Then $F$ is a Clairaut hemi-slant Riemannian map with $\tilde{s}=e^g$ if and only if
\begin{equation}\label{hemiNSC}
\begin{split}
 &g_N\Big(B(\nabla F_*)(X, {^{*}F_*BV})+\phi (F_*(\nabla^M_X {^{*}F_*BV}))-\phi(S_{CV}F_*X)\\&+B(\nabla^{F^\perp}_X CV)+B(\nabla^{F^\perp}_V CV)+\phi(\nabla^N_B CV), F_*X\Big) = -\frac{d(g \circ \beta)}{dt} g_N(F_*X, F_*X).   
\end{split}
\end{equation}
where $F_*X = F_*{X_1}+F_*{X_2},~F_*{X_1}\in\Gamma(D^{\perp}),~F_*{X_2}\in\Gamma(D^{\psi})$ and $V$ are the vertical and horizontal components of $\dot{\beta},$ respectively.
\end{theorem}
\begin{proof} Let $\beta$ be a geodesic on $N$ with constant speed $\sqrt{k},$ that is, $k = ||\dot{\beta}||^2.$
Let vertical and horizontal components of $\dot{\beta}$ be $F_*X$ and $V,$ respectively. Then, we get
\begin{equation}\label{hemisin}
g_N(F_*X, F_*X) = k\sin^2\psi(t) 
\end{equation}
and
\begin{equation*}
g_N(V, V) = k\cos^2\psi(t),
\end{equation*}
where $\psi(t)$ is the angle between $\dot{\beta}(t)$ and the horizontal subspace at $\beta(t).$ Differentiating \eqref{hemisin}, we get 
\begin{equation}\label{gNsincos}
g_N(\nabla^N_{\dot{\gamma}}F_*X, F_*X) = k \sin \psi(t) \cos \psi(t)\frac{d\psi}{dt}.
\end{equation}
Since $F$ is a hemi-slant Riemannian map, we can write $F_*X = F_*{X_1}+F_*{X_2},~F_*{X_1}\in\Gamma(D^{\perp}),~F_*{X_2}\in\Gamma(D^{\psi}).$ Then from the above equation, we have 
\begin{equation*}
g_N(\nabla^N_{\dot{\gamma}}F_*X, F_*X) = g_N(\nabla^N_{\dot{\gamma}}(F_*{X_1}+F_*{X_2}), F_*{X}).
\end{equation*}
From \eqref{gNsincos} and the above equation, we get
\begin{equation*}
k \sin \psi(t) \cos \psi(t)\frac{d\psi}{dt} = g_N(\nabla^N_{F_*X}F_*{X_1}+\nabla^N_{F_*X}F_*{X_2}+\nabla^N_V F_*{X_1}+\nabla^N_V F_*{X_2}, F_*X)
\end{equation*}
which gives
\begin{equation*}
\begin{split}
k \sin \psi(t) \cos \psi(t)\frac{d\psi}{dt} &= g_N(-J\nabla^N_{F_*X}JF_*{X_1}-J\nabla^N_{F_*X}JF_*{X_2}-J\nabla^N_V JF_*{X_1}\\&-J\nabla^N_V JF_*{X_2}, F_*X).
\end{split}
\end{equation*}
Making use of \eqref{Sv} and \eqref{hemiJV} in the above equation, we obtain
\begin{equation*}
\begin{split}
k \sin \psi(t) \cos \psi(t)\frac{d\psi}{dt} &= g_N\big(-J(-S_W F_*X+\nabla^{F^\perp}_X W)-J(\nabla^N_{F_*X}(\phi F_*{X_2}+\omega F_*{X_2}))\\&-J(\nabla^N_V W)-J(\nabla^N_V(\phi F_*{X_2}+\omega F_*{X_2})), F_*X\big),
\end{split}
\end{equation*}
where $JF_*{X_1} = W.$
Using \eqref{Sv}, \eqref{SFF}, \eqref{hemiJV} and \eqref{hemiJxi} in the above equation, we get
\begin{equation*}
\begin{split}
k \sin \psi(t) \cos \psi(t)\frac{d\psi}{dt} &= g_N\big(\phi(S_W F_*X)+\omega (S_W F_*X)-B(\nabla^{F^\perp}_X W)-C(\nabla^{F^\perp}_X W)\\&+\cos^2\theta \big(F_*(\nabla^M_X X_2)+(\nabla F_*)(X, X_2)\big)+S_{\omega(\phi F_*{X_2})}F_*X\\&-\nabla^{F^\perp}_X {\omega(\phi F_*{X_2})}+\phi(S_{\omega F_*{X_2}}F_*X)+\omega (S_{\omega F_*{X_2}}F_*X)\\&-B(\nabla^{F^\perp}_X \omega F_*{X_2})-C(\nabla^{F^\perp}_X \omega F_*{X_2})-B(\nabla^{F^\perp}_V W)\\&-C(\nabla^{F^\perp}_V W)+\cos^2\theta \nabla^N_V F_*{X_2}-\nabla^{F^\perp}_V \omega(\phi F_*{X_2})\\&-B(\nabla^{F^\perp}_V \omega F_*{X_2})-C(\nabla^{F^\perp}_V \omega F_*{X_2}), F_*X\big)
\end{split}
\end{equation*}
which implies
\begin{equation*}
\begin{split}
k \sin \psi(t) \cos \psi(t)\frac{d\psi}{dt} &= g_N\big(\phi(S_W F_*X)-B(\nabla^{F^\perp}_X W)+\cos^2\theta (F_*(\nabla^M_X X_2))\\&+S_{\omega(\phi F_*{X_2})}F_*X+\phi(S_{\omega F_*{X_2}}F_*X)-B(\nabla^{F^\perp}_X \omega F_*{X_2})\\&-B(\nabla^{F^\perp}_V W)+\cos^2\theta \nabla^N_V F_*{X_2}-B(\nabla^{F^\perp}_V \omega F_*{X_2}), F_*X\big).
\end{split}
\end{equation*}
Making use of Theorem \eqref{hemi-geodesic} in the above equation, we get
\begin{equation}\label{hemisincos}
\begin{split}
k \sin \psi(t) \cos \psi(t)\frac{d\psi}{dt} &= g_N\big(B(\nabla F_*)(X, {^{*}F_*BV})+\phi (F_*(\nabla^M_X {^{*}F_*BV}))-\phi(S_{CV}F_*X)\\&+B(\nabla^{F^\perp}_X CV)+B(\nabla^{F^\perp}_V CV)+\phi(\nabla^N_B CV), F_*X\big).
\end{split}
\end{equation}
Further, $F$ is a Clairaut hemi-slant Riemannian map with $\tilde{s} = e^g$ if and only if  
\begin{equation*}
\frac{d}{dt}( e^{g \circ \beta} \sin \psi(t)) = 0.
\end{equation*}
which implies  
\begin{equation*}
e^{g\circ\beta}\cos\psi(t)\frac{d\psi}{dt}+e^{g\circ\beta}\frac{d(g\circ \beta)}{dt}\sin\psi(t) = 0.
\end{equation*}
Since $e^{g \circ \beta}$ is a positive function, we get 
\begin{equation*}
\frac{d(g \circ \beta)}{dt} \sin\psi(t) + \cos\psi(t) \frac{d\psi}{dt} = 0. \end{equation*}
Multiplying the above equation with $k\sin\psi$ and comparing with \eqref{hemisincos}, we get \eqref{hemiNSC}.
\end{proof}
\begin{theorem}
Let $F:(M, g_M)\rightarrow (N, g_N, J)$ be a Clairaut hemi-slant Riemannian map with $\tilde{s} = e^g$ from a Riemannian manifold to a K\"ahler manifold such that \text{ker}$F_*$ is minimal. Then for $F_*X = F_*{X_1}+F_*{X_2}\in\Gamma(kerF_*)^\perp,$ where $F_*{X_1}\in\Gamma(D^{\perp})$ and $F_*{X_2}\in\Gamma(D^{\psi}),~F$ is a harmonic map if and only if
\begin{equation}
trace\Big(C(\nabla^{F\perp}_X W)+\nabla^{F\perp}_X\omega(\phi F_*{X_2})+C(\nabla^{F\perp}_X \omega F_*{X_2})\Big) = 0.
\end{equation}
\end{theorem}
\begin{proof}
We know that
\begin{equation*}
(\nabla F_*)(X, X) = \nabla^N_{F_{*}X}F_*X - F_{*}(\nabla^M_X X).
\end{equation*}
We can write $F_*X = F_*{X_1}+F_*{X_1},$ where $F_*{X_1}\in\Gamma(D^\perp)$ and $F_*{X_2}\in\Gamma(D^{\psi}).$ Making use of \eqref{nablaJ} in the above equation, we get
\begin{equation*}
\begin{split}
(\nabla F_*)(X, X) = -J\nabla^N_{F_{*}X}JF_*{X_1}-J\nabla^N_{F_{*}X}JF_*{X_2}- F_{*}(\nabla^M_X X).
\end{split}
\end{equation*}
Using  \eqref{Sv}, \eqref{SFF}, \eqref{hemiJV} and \eqref{hemiJxi} in above equation, we obtain
\begin{equation*}
\begin{split}
(\nabla F_*)(X, X) = &\phi(S_W F_*X)+\omega(S_W F_*X)-B(\nabla^{F^\perp}_X W)-C(\nabla^{F^\perp}_X W)+\cos^2\theta(\nabla F_*)(X, X_2)\\&+\cos^2\theta F_*(\nabla^M_X {X_2})+S_{\omega(\phi F_*{X_2})}F_*X-\nabla^{F^\perp}_X \omega(\phi F_*{X_2})+\phi(S_{\omega F_*{X_2}}F_*X)\\&+\omega(S_{\omega F_*{X_2}}F_*X)-B(\nabla^{F^\perp}_X \omega F_*{X_2})-C(\nabla^{F^\perp}_X \omega F_*{X_2})-F_*(\nabla^M_X X),
\end{split}
\end{equation*}
where $JF_*{X_1} = W.$ Since the range part of second fundamental form is zero, equating the $(range F_*)^\perp$ part of above equation and using Theorem \ref{NSC}, we get
\begin{equation*}
\begin{split}
-g_M(X, X)\nabla^N g =&-W(g)\omega F_*X-C(\nabla^{F^\perp}_X W)-\cos^2\theta g_M(X, X_2)\nabla^N g\\&-\nabla^{F^\perp}_X\omega(\phi F_*{X_2})-\omega((\omega F_*{X_2})(g)F_*X)-C(\nabla^{F^\perp}_X \omega F_*{X_2}).
\end{split}
\end{equation*}
Making use of Lemma \ref{harmonic} and taking trace in above equation, we get the required result.
\end{proof}
\begin{theorem}
Let $F:(M, g_M)\rightarrow (N, g_N, J)$ be a Clairaut hemi-slant Riemannian map with $\tilde{s} = e^g$ from a Riemannian manifold to a K\"ahler manifold. Then for any $F_*X, F_*Y\in\Gamma(D^\perp),$ $F_*Z\in\Gamma(D^{\psi})$ and $V\in\Gamma(rangeF_*)^\perp,$ the following assertions are equivalent:
\begin{itemize}
\item [(i)] The anti-invariant distribution $D^\perp$ is a totally geodesic foliation on $N,$
\item [(ii)]$JF_*Y(g)g_N(F_*X, \phi F_*Z)+g_N(\nabla^{F\perp}_X JF_*Y, \omega F_*Z) = 0$\\ and\\ $g_M(X, ^{*}F_* BV)g_N(\nabla^N g, JF_*Y)+g_N(\nabla^{F^\perp}_X JF_*Y, CV) = 0,$ 
\item [(iii)] $\omega(\phi F_*Z)(g)g_N(F_*X, F_*Y)+g_N(-g_M(X, Y)C(\nabla^N g)+\omega(F_*(\nabla^M_X Y)), \omega F_*Z) = 0.$ 
\end{itemize}
\end{theorem}
\begin{proof}
For any $F_*X, F_*Y\in\Gamma(D^\perp)$ and $F_*Z\in\Gamma(D^{\psi}),$ we have
\begin{equation*}
\begin{split}
 g_N(\nabla^N_{F_*X}F_*Y, F_*Z)& = g_N(\nabla^N_{F_*X}JF_*Y, JF_*Z).
 \end{split}
\end{equation*}
Making use of \eqref{Sv} and \eqref{hemiJV} in the above equation, we get
\begin{equation*}
g_N(\nabla^N_{F_*X}F_*Y, F_*Z) = g_N(-S_{J F_*Y}F_*X, \phi F_*Z)+g_N(\nabla^{F^\perp}_X JF_*Y, \omega F_*Z).
\end{equation*}
Using Theorem \ref{NSC} in the above equation, we get
\begin{equation}\label{semipsi}
\begin{split}
g_N(\nabla^N_{F_*X}F_*Y, F_*Z) = JF_*Y(g)g_N(F_*X, \phi F_*Z)+g_N(\nabla^{F^\perp}_X JF_*Y, \omega F_*Z).
\end{split}
\end{equation}
On the other hand, for $V\in\Gamma(range F_*)^\perp,$ we have 
\begin{equation*}
\begin{split}
g_N(\nabla^N_{F_*X}F_*Y, V) &=g_N(\nabla^N_{F_*X}JF_*Y, JV)\\& = g_N(\nabla^N_{F_*X}JF_*Y, BV)+ g_N(\nabla^N_{F_*X}JF_*Y, CV)\\& = -g_N(\nabla^N_{F_*X}BV, JF_*Y)+g_N(\nabla^{F^\perp}_{X}JF_*Y, CV)\\& = -g_N\big((\nabla F_*)(X, ^{*}F_* BV), JF_*Y\big)+g_N(\nabla^{F^\perp}_{X}JF_*Y, CV).
\end{split}
\end{equation*}
Making use of Theorem \ref{NSC} in the above equation, we get
\begin{equation}\label{semirperp}
g_N(\nabla^N_{F_*X}F_*Y, V) = g_M(X, ^{*}F_* BV)g_N(\nabla^N g, JF_*Y)+g_N(\nabla^{F^\perp}_X JF_*Y, CV).
\end{equation}
Now, from \eqref{semipsi}, \eqref{semirperp}, using Definition \ref{deftgf} the proof of (i) $\iff$ (ii) follows.\\
Next
\begin{equation*}
\begin{split} 
g_N(\nabla^N_{F_*X}F_*Y, F_*Z) &= g_N(\nabla^N_{F_*X}JF_*Y, JF_*Z)\\& = g_N(\nabla^N_{F_*X}JF_*Y, \phi F_*Z)+ g_N(\nabla^N_{F_*X}JF_*Y,\omega F_*Z)\\& = -g_N(\nabla^N_{F_*X}\phi F_*Z, JF_*Y)+ g_N(\nabla^N_{F_*X}JF_*Y,\omega F_*Z).
\end{split}
\end{equation*}
Making use of \eqref{Sv}, \eqref{hemiJV}, \eqref{hemiJxi} and Theorem \ref{hemiangle} in the above equation, we get
\begin{equation*}
\sin^2\theta g_N(\nabla^N_{F_*X}F_*Y, F_*Z) = -g_N(S_{\omega(\phi F_*Z)}F_*X, F_*Y)+g_N(C(\nabla F_*)(X, Y)+\omega(F_*(\nabla^M_X Y)), \omega F_*Z)
\end{equation*}
Using Theorem \ref{NSC} in above equation, we obtain
\begin{equation*}\label{semisin}
\begin{split}
\sin^2\theta g_N(\nabla^N_{F_*X}F_*Y, F_*Z) = &\omega(\phi F_*Z)(g)g_N(F_*X, F_*Y)+g_N(-g_M(X, Y)C(\nabla^N g)\\&+\omega(F_*(\nabla^M_X Y)), \omega F_*Z)  
\end{split}
\end{equation*}
From the above equation, using Definition \ref{deftgf}, the proof of $(i)\iff (iii)$ follows.
\end{proof}
\begin{theorem}
Let $F:(M, g_M)\rightarrow (N, g_N, J)$ be a Clairaut hemi-slant Riemannian map with $\tilde{s} = e^g$ from a Riemannian manifold to a K\"ahler manifold. Then, for any $F_*X, F_*Y\in\Gamma(D^{\psi})$ and $F_*Z\in\Gamma(D^\perp),$ the following assertions are equivalent:
\begin{itemize}
\item [(i)] The slant distributions $D^{\psi}$ is a totally geodesic foliation on $N.$
\item [(ii)] $$g_M(X, ^{*}F_*\phi F_*Y)g_N(\nabla^N g, \omega F_*Z) =g_N(\nabla^{F\perp}_X \omega F_*Y, \omega F_*Z)$$ and
 \begin{equation*}
\begin{split}
&g_N(\omega F_*(\nabla^M_X {^{*}F_*\phi F_*Y}), V)-g_M(X, {^{*}F_*\phi F_*Y})g_N(C(\nabla^N g), V) \\&= \omega F_*Y(g)g_N(F_*X, BV)+g_N(\nabla^{F\perp}_X \omega F_*Y, CV)
 \end{split}
 \end{equation*}
\item [(iii)] $$\omega F_*Y(g)g_N(\phi F_*X, F_*Z)=g_N\big(B(\nabla^{F\perp}_{F_*X}\omega(\phi F_*Y)), F_*Z\big).$$
\end{itemize}
\end{theorem}
\begin{proof}
Let $F:(M, g_M)\rightarrow (N, g_N, J)$ be a hemi-slant Riemannian map from a Riemannian manifold to a K\"ahler manifold. For any $F_*X, F_*Y\in\Gamma(D^{\psi})$ and $F_*Z\in\Gamma(D^\perp),$ we have
\begin{equation*}
\begin{split}
g_N(\nabla^N_{F_*X}F_*Y, F_*Z) &= g_N(\nabla^N_{F_*X}JF_*Y, JF_*Z)\\& = g_N(\nabla^N_{F_*X}(\phi F_*Y+\omega F_*Y), JF_*Z)\\& = g_N(\nabla^N_{F_*X} \phi F_*Y, JF_*Z)+ g_N(\nabla^N_{F_*X}\omega F_*Y, JF_*Z)\\& = g_N\big((\nabla F_*)(X, ^{*}F_*\phi F_*Y), JF_*Z)+g_N(-S_{\omega F_*Y}F_*X\\&+\nabla^{F^\perp}_X \omega F_*Y, JF_*Z).
\end{split}
\end{equation*}
Making use of Theorem \ref{NSC} in the above equation, we get
\begin{equation}\label{hemislantfoli}
g_N(\nabla^N_{F_*X}F_*Y, F_*Z) = -g_M(X, ^{*}F_*\phi F_*Y)g_N(\nabla^N g, \omega F_*Z)+g_N(\nabla^{F^\perp}_X \omega F_*Y, \omega F_*Z)
\end{equation}
On the other hand
\begin{equation*}
\begin{split}
g_N(\nabla^N_{F_*X}F_*Y, F_*Z) &= g_N(\nabla^N_{F_*X}JF_*Y, JF_*Z)\\& = g_N(\nabla^N_{F_*X}(\phi F_*Y+\omega F_*Y), JF_*Z)\\& = -g_N(J\nabla^N_{F_*X}\phi F_*Y, V)+g_N(\nabla^N_{F_*X}\omega F_*Y, JV)\\&= -g_N(C(\nabla F_*)(X, {^{*}F_*\phi F_*Y)+\omega F_*(\nabla^M_X{^{*}F_*\phi F_*Y}}), V)\\&-g_N(S_{\omega F_*Y}F_*X, BV)+g_N(\nabla^{F\perp}_X \omega F_*Y, CV).
\end{split}
\end{equation*}
Using Theorem \ref{NSC} in the above equation, we obtain
\begin{equation*}
\begin{split}
g_N(\nabla^N_{F_*X}F_*Y, F_*Z) &= g_M(X, {^{*}F_*\phi F_*Y})g_N(C(\nabla^N g), V)
-g_N(\omega F_*(\nabla^M_X {^{*}F_*\phi F_*Y}), V)\\&
+\omega F_*Y(g)g_N(F_*X, BV)+g_N(\nabla^{F\perp}_X \omega F_*Y, CV).
\end{split}
\end{equation*}
From \eqref{hemislantfoli}, the above equation and Definition \ref{deftgf}, we get $(i)\iff (ii).$\\
For $F_*X, F_*Y\in\Gamma(D^\perp),~F_*Z\in\Gamma(D^\psi),$ we have
\begin{equation*}
\begin{split}
g_N(\nabla^N_{F_*X}F_*Y, F_*Z) &= -g_N(J\nabla^N_{F_*X}JF_*Y, F_*Z)\\& = g_N(J\nabla^N_{F_*X}(\phi F_*Y+\omega F_*Y), F_*Z)\\& =- g_N(\nabla^N_{F_*X}(\phi^2F_*Y+\omega(\phi F_*Y)), F_*Z)+g_N(J\nabla^N_{F_*X}\omega F_*Y, F_*Z)\\&=\cos^2\theta g_N(\nabla^N_{F_*X}F_*Y, F_*Z)+g_N(S_{\omega(\phi F_*Y)}F_*X, F_*Z)\\&+g_N\big(J(-S_{\omega F_*Y}F_*X+\nabla^{F^\perp}_{X} \omega(\phi F_*Y)), F_*Z\big) 
\end{split}
\end{equation*}
With the help of above equation, \eqref{Sv}, \eqref{hemiJV}, \eqref{hemiJxi} and Theorem \ref{NSC}, we obtain
\begin{equation*}\label{otherhemislantfoli}
\sin^2\theta g_N(\nabla^N_{F_*X}F_*Y, F_*Z)=-(\omega F_*Y)(g)g_N(\phi F_*X, F_*Z)+g_N\big(B(\nabla^{F^\perp}_{X}\omega(\phi F_*Y)), F_*Z\big)
\end{equation*}
From the above equation and using Definition \ref{deftgf}, the proof of $(i)\iff(iii)$ follows.
\end{proof}
From above two theorems, we have the following decomposition Theorem.
\begin{theorem}
Let $F:(M, g_M)\rightarrow (N, g_N, J)$ be a Clairaut hemi-slant Riemannian map with $\tilde{s} = e^g$ from a Riemannian manifold to a K\"ahler manifold with integrable distribution $(range F_*).$ Then for $F_*X, F_*Y\in\Gamma(range F_*), F_*Z\in\Gamma(D^\psi)$ and $V\in\Gamma(range F_*)^\perp,$ the leaf of $(range F_*)$ is a locally product Riemannian manifold $M_{\perp} \times M_{\psi}$ if and only if
\begin{equation*}
(\omega(\phi F_*Z))(g)g_N(F_*X, F_*Y)+g_N(-g_M(X, Y)C(\nabla^N g)+\omega(F_*(\nabla^M_X Z)), \omega F_*Z) = 0
\end{equation*}
and
\begin{equation*}
\omega F_*Y(g)g_N(\phi F_*X, F_*Z)=g_N\big(B(\nabla^{F\perp}_{F_*X}\omega(\phi F_*Y)), F_*Z\big),  
\end{equation*}
where $M_\perp$ and $M_\psi$ denote the leaves of $D_\perp$ and $D_\psi,$ respectively
\end{theorem}
\begin{theorem}
Let $F:(M, g_M)\rightarrow (N, g_N, J)$ be a Clairaut hemi-slant Riemannian map with $\tilde{s} = e^g$ from a Riemannian manifold to a K\"ahler manifold. Then $F$ is totally geodesic if and only if the following conditions are satisfied:
\begin{itemize}
\item [(i)] $(kerF_*)$ and  $(kerF_*)^\perp$ are a totally geodesic distributions.
\item [(ii)] For $F_*X\in\Gamma(range F_*)$ and $F_*Y\in\Gamma(D^\psi),$ \begin{equation*}
\begin{split}
&\cos^2\theta\nabla^N_{F_*X}F_*Y-\omega(\phi F_*Y)(g)F_*X-\nabla^{F^\perp}_X \omega(\phi F_*Y)-\omega F_*Y(g)(\phi F_*X+\omega F_*X)\\&-B(\nabla^{F\perp}_X \omega F_*Y)-C(\nabla^{F\perp}_X \omega F_*Y)-F_*(\nabla^M_X Y) = 0.
\end{split}  
\end{equation*} 
\item [(iii)] for $F_*X\in\Gamma(range F_*)$ and $F_*Y\in\Gamma(D^\perp),$
$$W(g)(\phi F_*X+\omega F_*X)+B(\nabla^{F\perp}_X W)+C(\nabla^{F\perp}_X W)-F_*(\nabla^M_X Y) = 0.$$
\end{itemize}
\end{theorem}
\begin{proof}
A map $F$ is totally geodesic if and only if
\begin{equation*}
(\nabla F_*)(U, V) = 0,
\end{equation*}
\begin{equation*}
(\nabla F_*)(U, X) = 0,    
\end{equation*}
and
\begin{equation}\label{hemitgm}
 (\nabla F_*)(X, Y) = 0.   
\end{equation}
The first two equations directly imply (i). For the proof of part (ii) and (iii), we proceed as follows:\\
For $F_*X\in\Gamma(range F_*)$ and $F_*Y\in\Gamma(D^\psi),$ we have
\begin{equation*}
\begin{split}
(\nabla F_*)(X, Y) &= \nabla^N_{F_*X}F_*Y-F_*(\nabla^M_X Y)\\& =-J\nabla^N_{F_*X}JF_*Y-F_*(\nabla^M_X Y)\\& = -J\nabla^N_{F_*X}(\phi F_*Y+\omega F_*Y)-F_*(\nabla^M_X Y)\\& = -\nabla^N_{F_*X}(\phi^2F_*Y+\omega(\phi F_*Y))-J(\nabla^N_{F_*X}\omega F_*Y)-F_*(\nabla^M_X Y)\\& = cos^2\theta\nabla^N_{F_*X}F_*Y+S_{\omega(\phi F_*Y)}F_*X-\nabla^{F^\perp}_X \omega(\phi F_*Y)-J(-S_{\omega F_*Y}F_*X\\&+\nabla^{F^\perp}_X \omega F_*Y)-F_*(\nabla^M_X Y) 
\end{split}   
\end{equation*}
From \eqref{hemiJV}, \eqref{hemiJxi}, using Theorem \ref{NSC} in the above equation, we get
\begin{equation*}
\begin{split}
(\nabla F_*)(X, Y) &= \cos^2\theta\nabla^N_{F_*X}F_*Y-\omega(\phi F_*Y)(g)F_*X-\nabla^{F^\perp}_X \omega(\phi F_*Y)-\omega F_*Y(g)(\phi F_*X+\omega F_*X)\\&-B(\nabla^{F^\perp}_X \omega F_*Y)-C(\nabla^{F^\perp}_X \omega F_*Y)-F_*(\nabla^M_X Y).    
\end{split}
\end{equation*}
From \eqref{hemitgm} and the above equation, we get the proof of part (ii).\\
Next for $F_*X\in\Gamma(range F_*)$ and $F_*Y\in\Gamma(D^\perp),$
\begin{equation*}
\begin{split}
(\nabla F_*)(X, Y) &= \nabla^N_{F_*X}F_*Y-F_*(\nabla^M_X Y)\\&=J\nabla^N_{F_*X}JF_*Y-F_*(\nabla^M_X Y)\\&=J\nabla^N_{F_*X}W-F_*(\nabla^M_X Y)\\&=J(-S_W F_*X+\nabla^{F^\perp}_X W)-F_*(\nabla^M_X Y)
\end{split}
\end{equation*}
With the help of \eqref{hemiJV}, \eqref{hemiJxi} and  Theorem \ref{NSC}, we get
\begin{equation*}
\begin{split}
(\nabla F_*)(X, Y) = W(g)(\phi F_*X+\omega F_*X)+B(\nabla^{F^\perp}_X W)+C(\nabla^{F^\perp}_X W)-F_*(\nabla^M_X Y).
\end{split}
\end{equation*}
From \ref{hemitgm} and the above equation, we get the proof of part (iii).

\end{proof}
\begin{theorem}
Let $F:(M, g_M)\rightarrow (N, g_N, J)$ be a Clairaut hemi-slant Riemannian map with $\tilde{s} = e^g$ from a Riemannian manifold to a K\"ahler manifold. Then, for $F_*X, F_*Y\in\Gamma(D^{\psi}), F_*Z\in\Gamma (D^\perp),$ the distribution $D^{\psi}$ is integrable if and only if
\begin{equation}
g_N\big(\omega F_*Y(g)\phi F_*X+B(\nabla^{F\perp}_X \omega F_*Y), F_*Z\big) = g_N(\nabla^{F\perp}_Y \omega F_*X, JF_*Z).
\end{equation}
\end{theorem}
\begin{proof}
For $F_*X, F_*Y\in\Gamma(D^{\psi})$ and $F_*Z\in\Gamma (D^\perp),$ we have
\begin{equation*}
\begin{split}
g_N([F_*X, F_*Y], F_*Z) &= g_N(\nabla^N_{F_*X}F_*Y-\nabla^N_{F_*Y}F_*X, F_*Z)\\&=g_N(\nabla^N_{F_*X}JF_*Y, JF_*Z)-g_N(\nabla^N_{F_*Y}JF_*X, JF_*Z).
\end{split}
\end{equation*}
From \eqref{nablaJ}, \eqref{hemiJV} and the above equation, we obtain
\begin{equation*}
\begin{split}
g_N([F_*X, F_*Y], F_*Z) &= -g_N\big(\nabla^N_{F_*X}(\phi^2F_*Y+\omega(\phi F_*Y)), F_*Z\big)-g_N(J\nabla^N_{F_*X}\omega F_*Y, F_*Z)\\&+g_N\big(\nabla^N_{F_*Y}(\phi^2F_*X+\omega(\phi F_*X)), F_*Z\big)-g_N(\nabla^N_{F_*Y}\omega F_*X, JF_*Z).
\end{split}
\end{equation*}
Using \eqref{Sv}, \eqref{hemiJV} and Theorem \eqref{hemiangle} in the above equation, we get
\begin{equation*}
\begin{split}
g_N([F_*X, F_*Y], F_*Z) &= \cos^2\theta g_N(\nabla^N_{F_*X}F_*Y, F_*Z)+g_N(S_{\omega(\phi F_*Y)}F_*X, F_*Z)\\&-g_N(-\phi(S_{\omega F_*Y}F_*X)+B(\nabla^{F\perp}_X \omega F_*Y), F_*Z)\\&-\cos^2\theta g_N(\nabla^N_{F_*Y}F_*X, F_*Z)-g_N(S_{\omega(\phi F_*X)}F_*Y, F_*Z)\\&+g_N(\nabla^{F\perp}_Y \omega F_*X, JF_*Z).
\end{split}
\end{equation*}
Making use of Theorem \ref{NSC} in the above equation, we get
\begin{equation*}
\begin{split}
\sin^2\theta g_N([F_*X, F_*Y], F_*Z) = &-g_N\big(\omega F_*Y(g) \phi F_*X+B(\nabla^{F^\perp}_X \omega F_*Y), F_*Z\big)\\&+g_N(\nabla^{F^\perp}_Y \omega F_*X, JF_*Z).
\end{split}
\end{equation*}
We get the required result from Theorem \ref{Frobenius} and the above equation.
\end{proof}
\begin{example}
Let $M= N = \mathbb{R}^6$ be Euclidean spaces with Riemannian metrics $g_M$ and $g_N,$ respectively, defined by
\begin{equation*}
 g_M = dx_1^2+ dx_2^2+e^{2x_4}dx_3^2+e^{2x_4}dx_4^2+dx_5^2+dx_6^2. 
\end{equation*}
\begin{equation*}
g_N = dy_1^2+ dy_2^2+e^{2x_4}dy_3^2+e^{2x_4}dy_4^2+dy_5^2+dy_6^2.
\end{equation*}
Now, take the complex structure $J(y_1, y_2, y_3, y_4, y_5, y_6) = (-y_2, y_1, -y_4, y_3, -y_6, y_5).$
Then the bases of $T_pM$ is $\{e_1 = \frac{\partial}{\partial x_1}, e_2= \frac{\partial}{\partial x_2}, e_3 = e^{-x_4}\frac{\partial}{\partial x_3}, e_4 = e^{-x_4}\frac{\partial}{\partial x_4}, e_5 = \frac{\partial}{\partial x_5}, e_6 = \frac{\partial}{\partial x_6}\}$ and basis of $T_{F(p)}N$ is $\{e_1 = \frac{\partial}{\partial y_1}, e_2= \frac{\partial}{\partial y_2}, e_3 = e^{-x_4}\frac{\partial}{\partial y_3}, e_4 = e^{-x_4}\frac{\partial}{\partial y_4}, e_5 = \frac{\partial}{\partial y_5}, e_6 = \frac{\partial}{\partial y_6}\}.$\\
Define a map $F:(M, g_M)\rightarrow (N, g_N, J),$ given by
$$F(x_1, x_2, x_3, x_4, x_5, x_6) = (x_2\sin\alpha, 0,  x_3, 0, x_5, x_2\cos\alpha)$$
Then, we have
$$\text{ker}F_* = \text{span}\{e_1, e_4, e_6 \},$$
$$(\text{ker}F_*)^\perp = \text{span}\{X_1= e_2, X_2= e_3, X_3 = e_5\},$$
$$\text{range} F_* = \text{span}\{F_*{X_1}= \sin\alpha e_1'+\cos\alpha e_6', F_*{X_2} = e_3', F_*{X_3} = e_5'\},$$
$$(\text{range} F_*)^\perp = \text{span}\{V_1 = -\cos\alpha e_1'+\sin\alpha e_6', V_2 = e_2', V_3 = e_4'\}.$$
For $X\in\Gamma(\text{ker}F_*)^\perp,$ it can be seen that $g_M(X, X) = g_N(F_*X, F_*X)~~~$ which implies $F$ is a Riemannian map.
Moreover, we have $JF_*{X_1} = \sin\alpha e_2'-\cos\alpha e_5', JF_*{X_2} = e_4', JF_*{X_3} = e_6', J(V_1) = -\cos\alpha e_2'-\sin\alpha e_5', J(V_2) = -e_1', J(V_3) = -e_3'.$ Here $D^\perp = \text{span}\{JF_*{X_2}\}$ and $D^\psi = \text{span}\{JF_*{X_1},\\ JF_*{X_3}\}.$
Therefore $F$ is a hemi-slant Riemannian map with hemi-slant angle $\alpha.$
The non-zero Christoffel symbols are $\Gamma_{3 3}^{4} = -1, \Gamma_{4 4}^4 = 1, \Gamma_{3 4}^3 = 1 = \Gamma_{4 3}^3.$
For any $\Tilde{X} = a_1e_1'+ a_2e_2'+ a_3e_3'+ a_4e_4'+ a_5e_5'+ a_6e_6'$ and $\Tilde{Y} = b_1e_1'+ b_2e_2'+ b_3e_3'+ b_4e_4'+ b_5e_5'+ b_6e_6',$ we have $\nabla^N_{\Tilde{X}}J\Tilde{Y} = J\nabla^N_{\Tilde{X}}\Tilde{Y}$ which shows that the K\"ahlerian condition holds.\\
Further, for $X\in(\text{ker}F_*)^\perp,$
$$g_M(X, X) = g_M(a_1e_2+a_2e_3+a_3e_5, a_1e_2+a_2e_3+a_3e_5)= a_1^2+a_2^2+a_3^2.$$
and $$(\nabla F_*)(X, X) = -\mu_1e_1'+\mu_2(-\cos\alpha e_2'-\sin\alpha e_5')-\mu_3e_3'.$$ Therefore, $(\nabla F_*)(X, X) = -g_M(X, X)\nabla^N g$~for\\
$$g = \frac{-\mu_1y_1+\mu_2(-\cos\alpha y_2-\sin\alpha y_5)-\mu_3y_3}{a_1^2+a_2^2+a_3^2}.$$ This shows that $F$ is a Clairaut hemi-slant Riemannian map.
\end{example}
\section{Inequalities on semi-slant Riemannian map to K\"ahler manifolds}
In this section, we study some important inequalities (Casorati curvature inequality, Chen's first inequality) in the context of semi-slant and hemi-slant Riemannian maps.
\begin{lemma}\cite{cas cur}\label{lemcascur}
Let $F:(M, g_M)\rightarrow (N, g_N, J, c)$ be a Riemannian map from a Riemannian manifold to a complex space form with $\text{rank}$ $r\geq  3.$ Then the Casorati curvature inequality is given by
\begin{equation*}
\begin{split}
&\rho^{\mathcal{H}}\leq \delta^{\mathcal{H}}_C(r-1)+\frac{c}{4}+\frac{3c}{4r(r-1)}\sum_{i, j=1}^{r}\big(g_N(F_*{e_i}, JF_*{e_j})\big)^2,\\& \rho^{\mathcal{H}}\leq \hat{\delta}^{\mathcal{H}}_C(r-1)+\frac{c}{4}+\frac{3c}{4r(r-1)}\sum_{i, j=1}^{r}\big(g_N(F_*{e_i}, JF_*{e_j})\big)^2.
\end{split}
\end{equation*}
At a point $p\in M,$ the equality holds in any of the above two inequalities if and only if for suitable orthonormal bases, the following hold.
$$B^{\mathcal{H}^\alpha}_{11} = B^{\mathcal{H}^\alpha}_{22}=,\dots, = B^{\mathcal{H}^\alpha}_{r-1~ r-1} = \frac{1}{2}B^{\mathcal{H}^\alpha}_{rr},$$
$$B^{\mathcal{H}^\alpha}_{ij} = 0,~~~ 1\leq i\neq j\leq r.$$
\end{lemma}
\begin{theorem}
Let $F:(M, g_M)\rightarrow (N, g_N, J, c)$ be a semi-slant Riemannian map from a Riemannian manifold to a complex space form with \text{rank}~$r\geq 3.$ Then the Casorati curvature inequality takes the form
\begin{equation*}
\rho^{\mathcal{H}}\leq \delta^{\mathcal{H}}_C(r-1)+\frac{c}{4}+\frac{3c}{4r(r-1)}(2r_1+2r_2\cos^2\theta),
\end{equation*}
\begin{equation*}
\rho^{\mathcal{H}}\leq \hat{\delta}^{\mathcal{H}}_C(r-1)+\frac{c}{4}+\frac{3c}{4r(r-1)}(2r_1+2r_2\cos^2\theta).
\end{equation*}
The equality cases are the same as the equality cases of Lemma \ref{lemcascur}.
\end{theorem}
\begin{proof}
Let $F:(M, g_M)\rightarrow (N, g_N, J, c)$ be a semi-slant Riemannian map from a Riemannian manifold to a complex space form with \text{rank}~$r\geq 3.$ Now, we consider an orthonormal frame $\{F_*{e_1}, F_*{e_2} = JF_*{e_1},\dots, F_*{2r_1} = JF_* e_{2r_1-1}, F_*e_{2r_1+1},  F_*e_{2r_1+2} = \sec\theta\phi e_{2r_1+1}, \dots, F_*{e_r} = F_*e_{2r_1+2r_2-1}\}$ of $\Gamma(range F_*),$ where $\{F_*{e_1}, F_*{e_2} = JF_*{e_1},\dots, F_*{2r_1} = JF_* e_{2r_1-1}\}$ and $\{F_*e_{2r_1+1},\\  F_*e_{2r_1+2} = \sec\theta\phi e_{2r_1+1}, \dots, F_*{e_r} = F_*e_{2r_1+2r_2-1}\}$ are the orthonormal frame of $\Gamma(F_*D_1)$ and $\Gamma (F_*D_2),$ respectively.
Therefore, we obtain
\begin{equation*}
 \sum_{i, j=1}^{r}\big(g_N(F_*{e_i}, JF_*{e_j})\big)^2 = 2r_1+2r_2\cos^2\theta.   
\end{equation*}
With the help of Lemma \ref{lemcascur} and the above equation, we get the required result.
\end{proof}
\begin{lemma}\cite{chen1}\label{cheninq}
 Let $F:(M, g_M)\rightarrow (N, g_N, J, c)$ be a Riemannian map from a Riemannian manifold to a complex space form. Then the Chen's first inequality is
\begin{equation}
\begin{split}
K^{\mathcal{H}}{(\mathbb{P})} &\geq \frac{1}{2} \big\{2\rho^{\mathcal{H}}-\frac{(r-2)}{(r-1)}||\tau^{\mathcal{H}}||^2-c(r^2-r-2)-3c\big(\sum_{i, j=1}^{r}\big(g_N(F_*{e_i}, JF_*{e_j})\big)^2\\&-2(g_N(F_*{e_1}, JF_*{e_2}))^2\big)\big\},
\end{split}
\end{equation}   
where $\rho_{\mathcal{H}}$ is the scalar curvature defined on $(kerF_*)^\perp$ and equality holds if and only if there exists an orthonormal basis $\{e_i\}^{r}_{i=1}$ of $(kerF_*)^\perp$ and orthonormal basis $\{V_l\}^{2n}_{l=r+1}$ of $(rangeF_*)^\perp$ such that $B^{l}_{ij} = g_N((\nabla F_*)(e_i, e_j), V_l)$ and shape operator is 
\[
S_{r+1} =
\begin{pmatrix}
B_{11}^{r+1} & 0 & 0 & \cdots & 0 \\
0 & B_{22}^{r+1} & 0 & \cdots & 0 \\
0 & 0 & B_{11}^{r+1}+B_{22}^{r+1} & \cdots & \vdots \\
\vdots & \vdots & \vdots & \ddots & \vdots \\
0 & 0 & 0 & \cdots & B_{11}^{r+1}+B_{22}^{r+1}
\end{pmatrix}
\]

and

\[
S_{\ell} =
\begin{pmatrix}
B_{11}^{\ell} & B_{12}^{\ell} & 0 & \cdots & 0 \\
B_{12}^{\ell} & -B_{11}^{\ell} & 0 & \cdots & 0 \\
0 & 0 & \ddots & \cdots & 0 \\
\vdots & \vdots & \vdots & \ddots & \vdots \\
0 & 0 & 0 & \cdots & 0
\end{pmatrix}
\quad \text{for } \ell = r+2, \dots, 2n.
\]

\end{lemma}
\begin{theorem}
Let $F:(M, g_M)\rightarrow (N, g_N, J, c)$ be a semi-slant Riemannian map from a Riemannian manifold to a complex space form. Then the Chen's first inequality is
\begin{equation}
K^{\mathcal{H}}{(\mathbb{P})} \geq \frac{1}{2} \big\{2\rho^{\mathcal{H}}-\frac{(r-2)}{(r-1)}||\tau^{\mathcal{H}}||^2-c(r^2-r-2)-3c(2r_1+2r_2\cos^2\theta)-2\big\}
\end{equation}
\end{theorem}
The equality case is the same as the equality case of Lemma \ref{cheninq}.
\begin{proof}
Let $F:(M, g_M)\rightarrow (N, g_N, J, c)$ be a semi-slant Riemannian map from a Riemannian manifold to a complex space form with \text{rank}~$r\geq 3.$ Now, we consider an orthonormal frame $\{F_*{e_1}, F_*{e_2} = JF_*{e_1},\dots, F_*{2r_1} = JF_* e_{2r_1-1}, F_*e_{2r_1+1},  F_*e_{2r_1+2} = \sec\theta\phi e_{2r_1+1}, \dots, F_*{e_r} = F_*e_{2r_1+2r_2-1}\}$ of $\Gamma(range F_*),$ where $\{F_*{e_1}, F_*{e_2} = JF_*{e_1},\dots, F_*{2r_1} = JF_* e_{2r_1-1}\}$ and $\{F_*e_{2r_1+1},\\  F_*e_{2r_1+2} = \sec\theta\phi e_{2r_1+1}, \dots, F_*{e_r} = F_*e_{2r_1+2r_2-1}\}$ are the orthonormal frame of $\Gamma{(F_*D_1)}$ and $\Gamma{(F_*D_2)},$ respectively.\\
Therefore,
\begin{equation*}\sum_{i, j=1}^{r}\big(g_N(F_*{e_i}, JF_*{e_j})\big)^2 = 2r_1+2r_2\cos^2\theta, \end{equation*}  
\begin{equation*}
 (g_N(F_*{e_1}, JF_*{e_2}))^2 = 1.  
\end{equation*}
From above equations and Lemma \ref{cheninq}, we get required result.
\end{proof}
Proceeding in a similar manner, we obtain the following results.
\begin{theorem}
Let $F:(M, g_M)\rightarrow (N, g_N, J, c)$ be a hemi-slant Riemannian map from a Riemannian manifold to a complex space form with \text{rank}~$r\geq 3.$ Then the Casorati curvature inequality takes the form
\begin{equation*}
\rho^{\mathcal{H}}\leq \delta^{\mathcal{H}}_C(r-1)+\frac{c}{4}+\frac{3c}{4r(r-1)}(2r_2\cos^2\theta),
\end{equation*}
\begin{equation*}
\rho^{\mathcal{H}}\leq \hat{\delta}^{\mathcal{H}}_C(r-1)+\frac{c}{4}+\frac{3c}{4r(r-1)}(2r_2\cos^2\theta).
\end{equation*}
The equality cases are the same as the equality cases of Lemma \ref{lemcascur}.
\end{theorem}
\begin{theorem}
Let $F:(M, g_M)\rightarrow (N, g_N, J, c)$ be a hemi-slant Riemannian map from a Riemannian manifold to a complex space form. Then the Chen's first inequality takes the form
\begin{equation}
K^{\mathcal{H}}{(\mathbb{P})} \geq \frac{1}{2} \big\{2\rho^{\mathcal{H}}-\frac{(r-2)}{(r-1)}||\tau^{\mathcal{H}}||^2-c(r^2-r-2)-6cr_2\cos^2\theta)\big\}
\end{equation}
\end{theorem}
The equality case is the same as the equality case of Lemma \ref{cheninq}.
\section{Conclusion and Future scope}
In differential geometry and the theory of manifolds, geodesics are central to understand the geometric structure of spaces. In the context of manifolds, the notion of a Clairaut Riemannian map has been introduced as a generalization of Clairaut’s theorem. In this paper, We have extended this idea on semi-slant and hemi-slant Riemannian maps into a Kähler manifold, and studied geometric aspects of these maps and their associated distributions viz. conditions for total geodesicness and integrability. In future, one can further expand this notion by introducing analogues of Clairaut Riemannian maps for other geometric structures, such as contact or Sasakian manifold structures.
\section{Acknowledgement}
The first author is grateful for the financial support provided
by the CSIR (Council of Science and Industrial Research) Delhi, India. File
no.[09/1051(12062)/2021-EMR-I]. The corresponding author is thankful to the
Department of Science and Technology(DST) Government of India for providing
financial assistance in terms of FIST project(TPN-69301) vide the letter
with Ref No.:(SR/FST/MS-1/2021/104).\\

\end{document}